\documentclass[12pt]{amsart}
 \usepackage{amsmath}
\usepackage{amssymb,latexsym,cite}
\usepackage{amscd}
\usepackage{comment}
 \usepackage{xspace}
 \usepackage{verbatim}
 \usepackage[english]{babel}

 \newtheorem{theorem}{Theorem}[section]
 \newtheorem{lemma}[theorem]{Lemma}
 \newtheorem{proposition}[theorem]{Proposition}
 \newtheorem{definition}[theorem]{Definition}
 
 \newtheorem{corollary}[theorem]{Corollary}
 
 \newtheorem{remark}[theorem]{Remark}

     \allowdisplaybreaks

   \numberwithin{equation}{section}
     \numberwithin{theorem}{section}

    \title [Positive solutions of   semilinear equations and Schr\"odinger equations]{Positive solutions of a class of semilinear equations  with absorption and Schr\"odinger equations}
     \author{Alano Ancona}
\address{D\'{e}partement de math\'{e}matiques B\^at. 425
\\Universit\'{e} Paris Sud, Orsay (91405) France}
\email{alano.ancona@math.u-psud.fr}
\author{Moshe Marcus }
\address{Department of Mathematics, Technion\\
 Haifa 32000, ISRAEL}
 \email{marcusm@math.technion.ac.il}

    \date{{\small December 28, 2014}}

\begin{document}

          \maketitle
          
          \vspace{-5mm}
     	
\begin{center} D\'{e}partement de math\'{e}matiques B\^at. 425
\\ Universit\'{e} Paris Sud, Orsay (91405) France, \end{center}
\begin{center}Department of Mathematics, Technion\\
 Haifa 32000, ISRAEL. \end{center}

      \baselineskip=15truept
      
      \vspace{3mm}


      \vspace{2truemm}

       \renewcommand{\abstractname}

       \noindent{\bf ABSTRACT.} Several results about positive solutions -in a Lipschitz domain- of a nonlinear elliptic equation  in a general form $ \Delta  u(x)-g(x,u(x))=0$  are proved, extending thus some known facts in the case of  $ g(x,t)=t^q$, $q>1$, and a smooth domain. Our results include a characterization -in terms of a natural capacity- of a (conditional) removability property, a characterization of moderate solutions and of their boundary trace and a property relating arbitrary positive solutions to moderate solutions. The proofs combine techniques of non-linear p.d.e.\ with potential theoretic methods with respect to linear Schr\"odinger equations. A general result describing the measures that are diffuse with respect to certain capacities  is also established and used.
       The appendix by the first author provides classes of functions $g$  such that the nonnegative solutions of $ \Delta  u-g(.,u)=0$ have some ``good" properties that appear in the paper.

  \vspace{7truemm}

         \vskip 0.5\baselineskip

       \section{Introduction}

 In this paper we study positive solutions of semilinear elliptic equations with absorption of the form
\begin{equation}\label{geq}
- \Delta   u(x) + g(x, u(x))=0
\end{equation}
in a bounded Lipschitz domain $\Omega $ in ${ \mathbb  R}^N$, $N \geq 2$. A function $u\in L^1_{loc}(\Omega)$ is a solution of \eqref{geq} if $g\circ u\in L^1_{loc}(\Omega)$ and the equation holds in the distribution sense. Here
$$g\circ u(x):=g(x,u(x)) \quad \forall x\in \Omega.$$

\subsection{Assumptions} Our basic assumptions on $g$ are
\begin{equation}\label{g-basic}\begin{aligned}
&(i) &&g\in C(\Omega\times{\mathbb R}),\\
&(ii) &&\text{$g(x,\cdot)$ is odd and  increasing for every }   x\in \Omega, \\
&(iii) && \text{$g(x,\cdot)$  is convex on $[0,\infty)$} \;\forall x\in \Omega.
\end{aligned}\end{equation}

Notice that by (ii) and (iii) $g$ is differentiable with respect to $t$ at $t=0$. We put
$$h(x,t)=\begin{cases}g(x,t)/t &\text{if }t\neq 0\\ g'_t(x,0) &\text{if }t=0.\end{cases}$$
and in addition to (\ref{g-basic}) we assume:
\begin{quotation}\label{Vcondref}
$\exists \bar a>0$  such that, for every compact set $E\subset \partial \Omega $ and
every  positive solution $u$ of \eqref{geq} in $\Omega $ which is continuous in $\bar\Omega  \setminus  E$,
\begin{equation}\label{Vcond}
u=0    \quad   \text{on }\partial \Omega  \setminus  E   \Longrightarrow  h(x,u(x))\leq \bar a\, \mathrm{ dist} (x,E)^{-2}, \,\,  \forall  x\in \Omega .
\end{equation}

\end{quotation}

In particular, taking  $E=\partial \Omega $  this condition says that for every  positive solution $u$ of  \eqref{geq} in $\Omega $,
\begin{equation}\label{VcondSp}
  h(x,u(x))\leq \bar a\, \mathrm{ dist} (x,\partial \Omega )^{-2}, \,\,  \forall  x\in \Omega.
\end{equation}
As in \cite{MM-mod},  a basic ingredient in our study is the connection between (1.1) and the Schr\"odinger equation  \eqref{Sch}-\eqref{Vcond'} with potential $V$  given by  \eqref{Vdef}.  Condition  \eqref{VcondSp}  plays a  crucial role in this context.

 On using some set of assumptions on $g$, Appendix A describes concrete classes of absorption terms (not necessarily convex w.r.\ to  $t \in   (0,+ \infty )$) for which
(\ref{VcondSp}) or (\ref{Vcond}) holds. See Proposition \ref{ConcCond-St} and Theorem \ref{thap2bis}. It is also shown there   that  condition (\ref{Vcond}) follows from (\ref{VcondSp}) for some large classes of absorption terms $g$ (see Theorem \ref{thap2}).

Condition  \eqref{VcondSp} and Proposition \ref{ConcCond-St} are closely related to the Keller-Osserman
condition \cite{Kell,Oss}.  When $g$ \textit{is independent of the space
variable} the latter condition reads
\begin{equation}\label{KO}
  \psi(a)=\int_a^\infty \frac{dt}{\sqrt{2 G(t)}}<\infty, \quad \forall
a>0\quad {\rm where\ }G(t)=\int_0^tg(s)ds.
\end{equation}
 Under  condition  \eqref{KO} Keller and Osserman obtained a universal
estimate for solutions of  \eqref{geq} (see \cite[Theorem 1]{Kell}).  Numerous papers dealt with consequences of  the the Keller -- Osserman estimate and related results. Partial reviews and references can be found in \cite{Bandle,DDGR,Lm2011,OK+sub}.  Inequalities obtained  in \cite{DDGR} are used in Appendix A.

It was recently observed by one of us that, assuming \eqref{g-basic},  an alternative
proof of Proposition \ref{ConcCond-St}  can be
derived directly from an estimate of Keller \cite[inequality (25)]{Kell} which is, in fact, the basis for his main result. Moreover, in the case when
$g$ is independent of  $x$,
it can be shown that  assumption (A.2) also implies   \eqref{Vcond}
for general domains, not necessarily
Lipschitz. For details and related results see \cite[Section 5]{OK+sub}.

The basic example of an equation where our conditions are satisfied is
\begin{equation}\label{eq}
-  \Delta   u+|u|^{q-1}u=0,    \quad    q>1.
\end{equation}
It  is related to the study of  branching processes and superdiffusions (see Dynkin \cite{Dy91,Dbook1, Dbook2} and Le Gall \cite {LG93a,  LG93b, LG95, LGbook}) and has been intensively studied.

We mention two other interesting examples:
\begin{equation}\label{ex2}
(a) \quad g=\mathrm{dist}\,(\, \cdot \, , {\partial}\Omega)^\alpha g_0\,, \quad (b) \quad g=\psi_1^\alpha g_0 \quad \alpha\geq0,
\end{equation}
and $g_0$ satisfies  \eqref{g-basic} and (A.2), e.g.\
$g_0(t)=|t|^q\rm{sign}\, t$, $q \geq 1$ or $g_0(t)=(e^{|t|}-1)\rm{sign}\, t.$
Here $\psi_1$ denotes the first eigenfunction of $-\Delta$ in $\Omega$,  $\psi_1(x_0)=1$ at a point $x_0\in \Omega$. For a proof and  other  examples
see the Appendix.

\subsection {Connexion with Schr\"odinger equations} Our study of equation  \eqref{geq} employs non-linear elliptic p.d.e.\ techniques  combined with results on linear Schr\"odinger equations of the form
\begin{equation}\label{Sch}
  -  \Delta   u +Vu=0,
\end{equation}
where the potential $V\in L^ \infty _{loc}(\Omega )$ is nonnegative  and satisfies the condition
\begin{equation}\label{Vcond'}
   V(x)\leq \bar a \delta (x)^{-2},  \quad   \forall x\in \partial  \Omega,
\end{equation}
for some real $\overline a>0$. This class of potentials is denoted by $\mathcal{V}(\Omega,\bar a)$.
The point is that if $u$ is a positive solution of  \eqref{geq} then it satisfies  \eqref{Sch} with
\begin{equation}\label{Vdef}
   V(x)=h(x,u(x)) \quad \forall x\in \Omega
\end{equation}
and   \eqref{VcondSp} translates to  \eqref{Vcond'}.
We will write
$$L^V:=-\Delta+V.$$
A solution of equation  \eqref{Sch} is called an $L^V$ harmonic function. The term `harmonic' will be reserved for $L^V=-\Delta$.

Let $x_0\in \Omega   $ be a fixed reference point.  We will denote by $\mathbb K$ (resp.\ $\mathbb K^V$) the Martin kernel  with respect to  the operator $-  \Delta  $ (resp.\ $L^V:=- \Delta  +V$) in $\Omega $ normalized at $x_0$, i.e.\  $\mathbb K^V: \partial \Omega \ni  \zeta  \mapsto \mathbb K ^V_ \zeta  $ where $\mathbb K^V_ \zeta  $ is the unique positive $L^V$ harmonic function vanishing on $\partial \Omega \setminus  \{  \zeta   \} $  such that $ \mathbb K_ \zeta  ^V(x_0)=1$ (see Section 3 for more details).  Here we use property  \eqref{Vcond'}  that allows us  to apply the general results from \cite{Annals87}. If $\mu\in \mathfrak M(\partial \Omega )$ (the space of finite  Borel measures on $\partial \Omega $),  we denote by $ \mathbb K^V[\mu]$ (resp.\  $\mathbb K[\mu]$) the corresponding $L^V$-harmonic function (resp.\ $  \Delta  $-harmonic function).

It is known (ref.\ \cite{Annals87}) that every positive $L^V$ harmonic function $u$ in $\Omega$ can be uniquely represented in the form
\begin{equation}\label{Vrep}
  u={\mathbb K}^V[\mu]:=\int \mathbb K ^V_ \zeta  (.)\, d\mu  ( \zeta  )
\end{equation}
for some $\mu\in {\mathfrak M}_+(\partial\Omega)$. The measure $\mu$ will be called \emph{ the $L^V$ boundary measure} for $u$.

Let $ \mathbb G $ and $   \mathbb G ^V$ denote  the Green functions in  $\Omega$ of  $-  \Delta  $ and $L^V$ respectively.  If  $\mu  $ is a nonnegative Radon measure in $\Omega $ and  $x\in \Omega $, put
\begin{equation}\label{BBG}
 \mathbb G  [\mu](x):=\int_{\Omega   } \mathbb G (x,\xi )\,d\mu(\xi ),    \quad     \mathbb G  ^V[\mu](x):=\int_{\Omega } \mathbb G ^V(x,\xi )\,d\mu(\xi ).
\end{equation}

Note that if $ { \mathbb G }[\mu]$ is  finite at some point in $\Omega$ then it is locally integrable in $\Omega $ and this happens iff $\int_\Omega   \psi \, d\mu  < \infty $ where $\psi$ is the function $ \psi (x)=\mathbf 1\wedge \mathbb G (x_0,x)$, $x \in   \Omega $.  So this function plays a special role in a number of conditions below. One could use as well any continuous function $ \psi _1$ such that  $C^{-1}  \psi  \leq  \psi _1 \leq C \psi $ in  $ \Omega $ for some real $C \geq 1$.

We also note that if $ { \mathbb G }[\mu]\not \equiv +\infty $, then $ { \mathbb G }[\mu]$ is   superharmonic and has no positive  harmonic minorant in $\Omega $. Conversely a nonnegative superharmonic function $s$ in $\Omega $
that has no positive harmonic minorant is of the form $ s={ \mathbb G }[\mu ]$. For these  classical facts see e.g.\  \cite{brelotTCP}, \cite{AGpotT} or \cite{HLpoT}.

\subsection{Traces and moderate solutions of  \eqref{geq}} An exhaustion of $\Omega $ is any increasing sequence $\{\Omega _n\}$ of open and relatively compact subsets of $\Omega $ (so  $\bar\Omega_n\subset\Omega$) such that $\Omega = \cup _{n \geq 1} \Omega _n$. In what follows the reader may as well restrict   to exhaustions  by smooth domains. The next definition was introduced in \cite[Definition 3.6]{MVlip} w.r.\ to a  more restrictive concept of exhaustion.

\begin{definition} \label{mbtr} If $u$ is a nonnegative Borel function in $\Omega $, we say that $u$ has an $m$-boundary trace $\mu\in\mathfrak M_+(\partial \Omega )$ if, for any  exhaustion $\{\Omega   _n\}$ of $\Omega $,
\begin{equation}\label{mbtr1}
   \lim_{n\to\infty}\int_{\partial \Omega _n}u \varphi  \,d\omega _n=\int_{\partial \Omega } \varphi  \,d\mu,   \quad  \  \forall  \varphi  \in C(\overline   \Omega  ).
\end{equation}
where $\omega _n$ is $\omega _{\Omega _n}^{x_0}$, the harmonic measure of  $x_0$ in $\Omega _n$.
The m-boundary trace of $u$ on $\partial \Omega $ is denoted by $ {\mathrm tr}  _ {\partial \Omega } \,u$. (The subscript will be  omitted if it is clear from the context.)
\end{definition}

Recall that by definition (ref.\ \cite{brelot}, \cite{AGpotT}) the harmonic measure $\omega _{\Omega _n}^{x_0}$ is such that  for every $\varphi \in C(\partial \Omega _n)$, the integral $\int \varphi  \, d\omega _{\Omega _n}^{x_0}$ is the value  at $x_0$ of the solution $u$ to the Dirichlet problem $\Delta u=0$ in $\Omega _n$ and $u=\varphi $ on $\partial \Omega _n$.

Examples: (i)  If $u \in   C(\overline \Omega )$ then the $m$-boundary trace of $u_{ \vert \Omega }$ is the measure $u. \omega _\Omega ^{x_0}$ (since $\omega _n^{x_0}  \to  \omega _{\Omega }^{x_0}$ in the weak sense).
 {(ii)}   If  $u=\mathbb K[\nu]$, $\nu\in \mathfrak M_+(\partial \Omega )$ then $\nu$ is the m-boundary trace of $u$ (see  \cite[Lemma 2.2]{MVlip}).
 (iii) If $\mu\in \mathfrak M_{\psi}^+(\Omega )$ (i.e.\ $\int  \psi \, d\mu  < \infty $) then $\mathbb G  [\mu]$ has m-boundary trace zero. (See \cite[Lemma 3.1]{MVlip}).

 It follows from the Riesz decomposition theorem that every nonnegative superharmonic function in $\Omega $ admits a $m$-boundary trace which is the  $m$-boundary trace of its largest harmonic minorant. In general this statement does not apply to $L^V$ harmonic functions and  if $u=\mathbb K^V[\nu]$, $\nu  \in   \mathfrak M_+(\partial \Omega )$ then, in general, $\nu $ is \emph{not} the m-boundary trace of $u$.

\medskip

We will consider boundary value problems of the form
\begin{equation}\label{g-bvp}   \begin{aligned}
  -  \Delta   u+g\circ u&=0 &&\text{ in $\Omega   $},\\
   {\mathrm tr}  \,u&=\nu &&\text{ on $\partial \Omega $},
   \end{aligned}\end{equation}
where $\nu\in \mathfrak M_+(\partial \Omega )$ and $ {\mathrm tr}  \,u$ denotes the m-boundary trace of $u$.

The next definition provides a natural class of solutions of (\ref{geq}) related to such   problems. Recall that $\psi (x)=1\wedge G(x_0,x)$.
\begin{definition}\label{gmod-sol} Let $u$ be a nonnegative solution of (\ref{geq}). We say that $u$ is a $g$ moderate solution  of  (\ref{geq}) if
$g\circ u\in L^1_\psi(\Omega)$ (i.e.\, $\int _\Omega  g(x,u(x))\,  \psi (x)\, dx < \infty $).
\end{definition}

Notice that the condition above means that the Green potential $\mathbb G(-\Delta u)$ is not identically $+\infty $ in $\Omega $  and implies that $\mathbb G(-\Delta u)$ has boundary trace zero. See examples (iii) above.
 \begin{comment}
 Equivalently,  $\mu$ being the nonnegative measure $\Delta u$, there is a positive $L^1_{\rm loc}(\Omega )$ solution $s$ to $-\Delta s=\mu$, or -on taking  $s=\mathbb G(\mu )$- there is such an $s$ with moreover ${\mathrm tr}  \,u=0$.
 \end{comment}
The next proposition follows easily. It recalls among other things some equivalent definitions and the connection with the previous boundary value problems.

\begin{proposition}\label{gmodprop} If $u$ is a nonnegative solution of (\ref{geq}) then $u$ is $g$-moderate iff $u$ admits a harmonic majorant in $\Omega $. This condition is also equivalent to the existence of an $m$-boundary trace for $u$ on $\partial \Omega $. If $u$ is $g$ moderate then $u \in   L^1_ \psi (\Omega )$ and if $\nu $ is its $m$-boundary trace, $ { \mathbb  K}[\nu]$ is the least harmonic majorant of $u$ in $\Omega $.
\end{proposition}
There are also well-known equivalent variational definitions of solutions of  (\ref{g-bvp}).

\begin{remark}\label{rksuper}  If $u$ is a positive supersolution  of  \eqref{geq}  such that  $g\circ u\in L^1_\psi(\Omega )$ then $u$ has an m-boundary trace and the largest solution of  \eqref{geq} dominated by $u$ has the same m-boundary trace. There is a parallel statement for subsolutions dominated by a moderate solution. See \cite[Theorem 4.3]{MVlip}.
\end{remark}

It is known that, for every $\nu\in {\mathfrak M}_+(\partial\Omega)$, problem  \eqref{g-bvp} has at most one solution.
A measure $\nu$ for which a solution exists is called a \emph{$g$-good measure} and the corresponding solution of  \eqref{g-bvp} is denoted by ${\mathbb S}^g[\nu]$.
Moreover, if $\nu_1$, $\nu_2$ are $g$-good measures then,
\begin{equation}\label{monotone}
   \nu_1\leq \nu_2 \Longrightarrow {\mathbb S}^g[\nu_1]\leq {\mathbb S}^g[\nu_2].
\end{equation}
For a proof see e.g. \cite{MVlip} (where $g$ is assumed to be space independent).

We next recall some basic stability properties of good measures.

\begin{proposition}\label{gmodstab} (i) If $\mu,\, \nu    \in   \mathfrak M _+(\partial \Omega )$,  $\mu  $ is a $g$-good measure and $0 \leq \nu  \leq \mu  $ then $\nu $ is also $g$-good.  (ii) A convex combination of $g$-good measures is again $g$-good. (iii) If $M \subset   \mathfrak M_+(\partial \Omega )$ is a set of $g$-good measures dominated in $ \mathfrak M_+(\partial \Omega )$ then the least upper bound of $M$ is again $g$-good. \end{proposition}

We sketch some proofs for the reader's convenience. Under the assumptions in (i),  if $u= {\mathbb S}^g[\mu]$,  then $v=(u-{\mathbb K}[\mu -\nu])_+$ is a subsolution of (\ref{geq}) and $\mathbb K[\nu]$ is a supersolution larger than $v$.  So if $w$ denote the largest subsolution (\ref{geq})
dominated by ${\mathbb K}[\nu ] $, $w$ is  a solution and (i) follows from $ {\mathrm tr}  _ {\partial \Omega } (u-{ \mathbb K}[\mu -\nu  ])=\nu $, ${\mathrm tr}  _ {\partial \Omega } ({ \mathbb K}[\nu  ])=\nu$ and $u-{ \mathbb K}[\mu -\nu  ]\leq w\leq { \mathbb K}[\nu  ]$.
To prove (iii) 
 observe that  if $\lambda $ is the least upper bound of $M$, then $ {\mathbb K}[\lambda ]$ is the least harmonic majorant of $w=\sup _{\mu \in M} {\mathbb S}^g[\mu]$. If $u$ is  the smallest supersolution dominating $w$,   one may then either observe that  $p={ \mathbb K}[\lambda ]-u$ has no positive harmonic minorant and so $\mathrm{tr\,}(p)=0$ or note that  $u$ is a moderate solution (since $u\leq { \mathbb K}[\lambda ]) $ to conclude that $ {\mathrm tr}  _ {\partial \Omega } (u)=\lambda$.
Finally statement (ii) can be proved as in
\cite[Corollary 4.5]{BMP04}). $\square$ 

In general the sum of two $g$-good measures  need not be $g$-good. For instance this is the case when $g(x,t)=\sinh (t)$. But if $g$ satisfies the (uniform) $ \Delta  _2$ condition  --i.e.\  $g(x,2t) \leq C\, g(x,t)$, $ \forall  ( x ,t)\in   \Omega \times { \mathbb  R}_+$ for some $C \geq 0$-- then the sum of two $g$-good measures is $g$-good and the set of $g$-good measures is a convex cone.

Denote
\begin{equation}\label{CTg0}
  \mathfrak M_+^{g,  \Delta  }(\partial \Omega ):=\{\nu\in \mathfrak M_+(\partial \Omega ): \int_\Omega    (g\circ\mathbb K[\nu])\psi\,dx <\infty\}.
 \end{equation}
 Notice that this is a convex set  of good measures. If $\nu\in\mathfrak M_+^{g,  \Delta  }(\partial \Omega )$ we say that \emph{$\nu$ is a $(g,\Delta)$ - good measure}. When $g$ satisfies the $ \Delta  _2$-condition the set  $\mathfrak M_+^{g,  \Delta  }(\partial \Omega ) $  is a convex \emph{cone} of measures.

\subsection{Main results}
The first set of results concerns $g$-moderate solutions. In all of these we assume that $g$ satisfies our basic conditions, namely,  \eqref{g-basic} and  \eqref{Vcond}. In some of the results we assume in addition that $g$ satisfies the $\Delta_2$ condition.

For the statement of the results it is convenient to introduce the following definitions.

\begin{definition}\label{removable} A compact set $F \subset  \partial \Omega $ is \emph{conditionally $g$-removable} if, for any non-negative $g$-moderate solution $u$ of  \eqref{geq},
$$u\in C(\overline   \Omega  \setminus  F)\;\text{and } u=0\;\text{on }\partial\Omega\setminus F \;\Longrightarrow \;u=0.$$
\end{definition}

When $\mathcal T   \subset   { \mathfrak M}_+(\partial \Omega )$, we  say that a Borel set $A \subset  \partial \Omega $ is \emph{$\mathcal T  $-null}  if
\begin{equation}\label{T-null}
  \tau(A)=0\quad  \forall \tau\in \mathcal T.
\end{equation}

\begin{theorem}\label{I-1}
 A compact set $F  \subset    \partial \Omega $ is  conditionally $g$-removable  if and only if  $F$ is ${\mathfrak M}^{g,\Delta}_+$-null.
\end{theorem}

\begin{theorem}\label{I-2}
Let  $u$ be a positive $g$-moderate solution of  \eqref{geq} in $\Omega $. Then  there exists a sequence  $\{\mu_n\}\subset{\mathfrak M}^{g,\Delta}_+:={\mathfrak M}^{g,\Delta}_+(\partial \Omega )$  such that
\begin{equation}\label{I-2.1}
  u=\lim_{n\to\infty} {\mathbb S}^g[\mu_1+\cdots +\mu_n].
\end{equation}

If $g$ satisfies the $\Delta_2$ condition,  there is an increasing sequence  $\{\mu_n\}\subset{\mathfrak M}^{g,\Delta}_+$  such that
\begin{equation}\label{I-2.2}
    u=\lim_{n\to\infty}{\mathbb S}^g[\mu_n].
\end{equation}
\end{theorem}

\noindent\textit{Remark.}\hskip 2mm Note that in the first part of the theorem, $\mu_1+\cdots +\mu_n$ needs not be  $(g,\Delta)$-good. But, as it is dominated by ${\mathrm tr}(u)$, the sum is $g$-good.

\begin{theorem}\label{I-3}
(i) If $\mu\in    {\mathfrak M}_+:={\mathfrak M}_+(\partial \Omega )$ is a $g$-good measure then $\mu$ vanishes on every ${\mathfrak M}^{g,\Delta}_+$-null set.

\noindent(ii) If $\mu\in    {\mathfrak M}_+$ vanishes on ${\mathfrak M}^{g,\Delta}_+$ null sets then there exists a  Borel function $f:\partial\Omega\mapsto(0,1]$  such that $f\mu$ is a $g$-good measure.

\noindent(iii) If, in addition to the basic conditions, $g$ satisfies the $\Delta_2$ condition then, $\mu\in    {\mathfrak M}_+$ is a $g$-good measure if and only if  $\mu$ vanishes on every ${\mathfrak M}^{g,\Delta}_+$-null set.
\end{theorem}

Notice that $(ii)$ means that the greatest $g$-good measure $\nu $ smaller than $ \mu  $ is in the form $\nu =f\mu  $ with $f>0$ $\mu  $-a.e.

Theorems \ref{I-1}--\ref{I-3} are proved in section 5.

Much of the research on the problems discussed above focused on the case of power nonlinearities, i.e.\
equation  \eqref{eq} and the corresponding non-homogenous equation
\begin{equation}\label{eq'}
   -\Delta u + |u|^{q-1}u=\mu
\end{equation}
where $\mu\in {\mathfrak M}(\Omega)$.

Benilan and Brezis studied the Dirichlet problem for \eqref{eq'} and proved: (a) If $q\geq N/(N-2)$ and $\mu=\delta_{x_0}$ ($=$ Dirac measure at a point $x_0\in  {\Omega} $) then the problem has no solution.(b) If $1<q<\geq N/(N-2)$ (the, so called, subcritical case) the boundary value problem possesses a solution for every finite measure $\mu$. The results,  obtained in the 70's,  have been published  in full only  in \cite{BenBr}. In part they appeared in \cite{Br-note} and \cite{Br80}. Brezis and Veron \cite{BrV} showed that if $q\geq N/(N-2)$ then, for positive solutions, any isolated singularity in $\Omega$ is removable. Baras and Pierre \cite{BP84} provided a complete characterization of good measures and removable subsets of $\Omega$ when $q\geq N/(N-2)$. (Uniqueness holds under very general conditions, see \cite{Br-note}.)

Following these results, the research turned to boundary value problems for \eqref{eq}, when $\Omega$ is a domain of class $C^2$ and the  boundary data is a positive measure. Assuming that $ {\Omega} $ is a bounded domain of class $C^2$,
 Gmira and Veron \cite{GV} proved that, if
$1<q<q_c=(N+1)/(N-1)$, every measure in ${\mathfrak M}(\partial\Omega)$ is $g$-good and, if $q\geq q_c$, isolated singularities on the boundary are removable.

Regarding \eqref{eq} in the supercritical case, $q\geq q_c$, the results stated in Theorems \ref{I-1} -- \ref{I-3}  have been previously established,   for domains of class $C^2$,  in a series of papers:
Le Gall \cite{LG93b,LG95} proved these results in the case $q=2$,
 Dynkin and Kuznetsov \cite{DK96, DK98a, DK98b} in the case $q_c\leq q\leq 2$ (employing mainly  probability methods) and Marcus and Veron \cite{MVsuper} in the case $q\geq 2$ (using purely analytic methods). A unified approach that applies to all $q\geq q_c$ was provided in \cite{MVrem}.

In the case of Lipschitz domains, equation  \eqref{eq}  has been investigated in \cite{MVlip} under conditions implying that $q$ is subcritical. The supercritical case was treated in \cite{MVLip} in the case that the domain is a polyhedron.
For general Lipschitz domains, Theorems \ref{I-1} -- \ref{I-3} are new even in the case of power nonlinearities.

We note here that equation  \eqref{geq} with $g(x,t)=\delta(x)^\alpha|t|^{q-1}t$, $\alpha>-2$ has been studied in  \cite{MVcpam}, for domains of class $C^2$ and $q$ in the subcritical range,  $1<q<(N+\alpha+1)/(N-1)$. Equation  \eqref{geq} with $g(t)=e^t-1$ was studied in  \cite{LV_MZ} for domains of class $C^2$.

For the statement of our next result on arbitrary positive solutions of  \eqref{geq} (not necessarily $g$-moderate) we introduce some additional notation and terminology.

Denote by $C_{g,\Delta}$  the set function on the Borel field ${ \mathcal B}or( \partial \Omega )$ defined by
\begin{equation}\label{g-cap}
    C_{g,\Delta  }(F)=\sup\{\tau(F):\tau\in {\mathfrak M}_+( \partial \Omega ),\; \; \int_\Omega  \,(g\circ {  \mathbb  K} _\tau)\psi\,dx\leq 1\}
\end{equation} for every Borel set $F \subset  \partial \Omega $. This set function is a capacity (in the sense of \cite{Moko}).

The family of $\mathcal M_+^{g,  \Delta  }$ - null sets can also be described as the family of Borel sets $A\subset\partial\Omega$ such that $C_{g,\Delta}(A)=0$.

 A   measure $  \mu \in {\mathfrak M}_+( \partial \Omega )$ is said to be \emph{diffuse relative to $C_{g, \Delta  }$} if $\mu$ vanishes on $C_{g, \Delta  }$ null sets.
 The measure $\mu$ is \emph{concentrated (or singular) w.r.\ to $ C_{g, \Delta  }$} if it is concentrated on a $C_{g, \Delta  }$ null set.
It is known that every measure $\mu\in {\mathfrak M}_+(\partial\Omega)$ can be uniquely written  in the form $\mu=\mu_d+\mu_c$ where $\mu_d$ is diffuse and $\mu_c$ is singular w.r.\ to $C_{g,\Delta}$ (see Section 6).
\begin{theorem}  \label{I-gsmod}  Assume that $g$ satisfies  \eqref{g-basic} and  \eqref{Vcond}. Let $u$ be a positive solution of  \eqref{geq}, let $V=h\circ u$ and let $\mu\in {\mathfrak M}_+(\partial\Omega)$ be the $L^V$ boundary measure of $u$, i.e.\  $u={\mathbb K}^V[\mu]$.

If $\mu$ is not singular  relative to $C_{g,\Delta}$ then $u$ dominates a positive $g$ moderate solution $v$ such that $\mu_d$ is absolutely continuous w.r. to $\mathrm{tr}\,v$.
\end{theorem}
For the proof of this theorem and some additional results see Section 6.

 In connexion with the  characterization  of positive solutions of  \eqref{eq} in terms of their boundary behavior,  Dynkin and Kuznetsov \cite{DK98b} (see also \cite{Dbook1}) raised a central question:
is every positive solution of  \eqref{eq} $\sigma$-moderate, i.e.\ the limit of an increasing sequence of $g$-moderate solutions ?

The question was answered affirmatively  in the case of $C^2$ domains in a series of works:  first, in the subcritical case (Marcus and Veron \cite{MVcras96, MVsub}), then in the supercritical case,
 by Mselati \cite{Ms} for $q=2$ and Dynkin \cite{Dbook2} for $q_c\leq q\leq 2$ and finally by Marcus
 \cite{MM-mod} for every $q\geq q_c$.
A crucial step in the proof of \cite{MM-mod} was to show that if $u$ is a positive solution of  \eqref{eq} and $\mu$ is defined as in Theorem \ref{I-gsmod}, then $\mu$ is diffuse relative to $C_{g,\Delta}$ and $u$ dominates a positive $g$-moderate solution. Theorem \ref{I-gsmod} constitutes a first step in an investigation of possible extensions of \cite{MM-mod} to Lipschitz domains and to more general non-linearities.

 Our study relies on results of
 \cite{MM-mod}, here adapted to Lipschitz domains, and potential theoretic results about linear Schr\"odin\-ger equations, in particular a recent result of  \cite{An-MZ} and a boundary Harnack principle based on  \cite{Annals87}. See also Sections 3 and  4 where some auxiliary  results, needed in the later sections, are  proved.

We also use and prove a general result (see Section 2, Theorem 2.1) that characterizes the  positive Borel measures which are diffuse with respect to a capacity defined as the supremum of a family of measures.  It extends standard results, in particular results  from \cite{FdlP}, \cite{DMaso} and \cite{BP84}.

{\sl Acknowledgement.} The authors thank the referees for several remarks that helped to improve the exposition.

\section{A measure theoretic approximation result.}

Let $X$ be a metric space and let $\mathfrak M_+(X)$ denote the set of all  \emph{finite} nonnegative Borel measures on $X$.

Given $\mathcal T   \subset   { \mathfrak M}_+(X)$, we say that a Borel set $A \subset  X$ is $\mathcal T  $-null  if
\begin{equation}\label{T-removable}
  \tau(A)=0\quad  \forall \tau\in \mathcal T  .
\end{equation}

 \begin{theorem}\label{rich}  If $X$ is compact and    $\mathcal K  \subset   \mathfrak M_+(X)$ satisfies the conditions
\begin{equation}\label{rich-space}   \begin{aligned}
  (a)\; &  \mathcal K {\rm\ is\  nonempty, convex \ and\  weakly\  closed \ in\ }  { \mathfrak M}_+(X),  \\
  (b)\; &\nu\in \mathcal K,\;\tau\in \mathfrak M_+(X),\;\tau\leq \nu\; \Longrightarrow \; \tau\in \mathcal K,\\
   \end{aligned}\end{equation}
then every  $\mu\in \mathfrak M_+(X)$ that vanishes on $\mathcal K$-null sets is in the form $\mu  =\sup _{ n \geq 1}n\,\nu  _n$ where $\nu_n \in   { \mathcal K}$ for all $n \geq 1$ and $ \{ n\nu_n \}$ is an increasing sequence  in $\mathfrak M_+(X)$.

\end{theorem}

Note that $ { \mathcal E}_+:=\cup _{n \geq 1} n\mathcal K$ is the convex cone generated by ${ \mathcal K}$ and that the conclusion says that $\mu  =\lim_{n \to   \infty }  \lambda _n$ for some increasing sequence $ \{  \lambda _n \}$  in ${ \mathcal E}_+$.

This theorem extends well-known results, in particular results of Feyel--de La Pradelle \cite{FdlP} and Dal Maso \cite{DMaso} in which ${\mathcal K} $ generates the positive cone of the dual of a Dirichlet space (resp.\ a Sobolev space $W^{-   \alpha ,q}(\partial \Omega )$). In  the proof below we use a result due to G.\,Mokobodzki \cite{Moko}.
\medskip

To prove Theorem \ref{rich} we first notice that upon replacing ${ \mathcal K}$ by $ \tilde { \mathcal K}:= \{ \mu   \in   { \mathcal K}\,;\, \mu  (X) \leq 1\, \}$, we may assume that ${ \mathcal K}$ is \emph {weakly compact} in ${ \mathcal M}_+(X)$. Without loss of generality we  will also assume that $ \bigcup    _{\mu   \in   {\mathcal K} }  \mathrm { supp }(\mu  )$ is dense in $X$ (this assumption ensures that the function $\Vert.\Vert $ as defined below is actually a norm).

Define then a set function as follows. For any Borel   subset $A$ of $X$, set
$${ \mathrm {Cap}} (A):=\sup  \{  \int _A \, d\mu  \,;\,  \, \mu   \in   { \mathcal K}  \}.$$

${ \mathrm {Cap}}$ is a capacity on $X$ in the sense of  \cite{Moko}. That is ${ \mathrm {Cap}}$ is finite nonnegative  and
\begin{itemize}

\item[(i)] ${ \mathrm {{Cap}}}( \emptyset  )=0$,

\item[(ii)] ${ \mathrm {Cap}}( \bigcup    _{n \geq 1} A_n) \leq  \sum_{n \geq 1} { \mathrm {Cap}}(A_n)$ for every sequence $ \{ A_n \}$ of Borel subsets of $X$

\item[(iii)] ${ \mathrm {Cap}}( \bigcup    _{n \geq 1} A_n) = \sup _{n \geq 1} { \mathrm {Cap}}(A_n)$ for every increasing sequence $ \{ A_n \}$ of Borel subsets of $X$,

\end{itemize}

as easily verified.   So the result in \cite{Moko} says the following.

\begin{lemma} Let $\nu  \in   { \mathfrak M}_+(X)$ be such that $\nu (A)=0$ for every Borel $ { \mathcal K}$-null set. Then there exists an increasing sequence $ \{ \nu _n \}$ in ${ \mathfrak M}_+(X)$  such that $\nu =\sup_n \nu _n$ and  $\nu _n \leq n\, { \mathrm{Cap}}$ for all $n \geq 1$.
\end{lemma}

 Thus to prove Theorem \ref{rich} it suffices to show that every measure $ \nu  \in   { \mathfrak M}_+(X)$ majorized by $\mathrm  {Cap}$ (i.e.\ $\nu (A) \leq {\mathrm {Cap}}(A)$ for every Borel set $A \subset   X$) is the upper enveloppe of a set of  measures in ${ \mathcal E}_+={ \mathbb  R}_+. {\mathcal K} $.

To begin with, we consider ${ \mathcal E}={ \mathcal E}_+-{ \mathcal E}_+$  the vector subspace of ${ \mathfrak M}(X;{ \mathbb  R})$ generated by ${\mathcal K} $  and define below a topology on $E:={ \mathcal C}(X;{ \mathbb  R})$  that agrees with the duality $ (E;{ \mathcal E})$.

Set for each $f \in   E$:
$$ \Vert f \Vert =\sup  \{  \,\vert \! \int f\, d\mu \,\vert    \,;\,  \mu  \in {\mathcal K}-{\mathcal K}  \, \}. $$

This defines a norm $ \Vert . \Vert $ on $E$ and this norm is equivalent to the norm
$$ \vert   \vert   \vert   f  \vert   \vert   \vert  =\sup  \{  \int  \vert  f \vert  \, d\mu  \,;\, \mu   \in   {\mathcal K} \,  \},\ \ f \in   E.$$

In effect, $\mu   \in  {\mathcal K}-{\mathcal K} $ if and only if $\mu  ^+$, $\mu  ^- \in    {\mathcal K}$ (using  (b) in
  \eqref{rich-space} and $\mu  =\mu  ^+-\mu  ^-$). Thus for $f \in  E$, $\mu   \in   { \mathcal K}-{ \mathcal K}$, we have  $ \vert  \int f\, d\mu   \vert  \leq \int  \vert  f \vert   \,d\mu  ^+ +\int  \vert  f \vert   \,d\mu  ^-$ and  $ \Vert f \Vert  \leq  2\, \vert   \vert   \vert  f \vert   \vert   \vert  $. On the other hand if  $\mu   \in   { \mathcal K} $ and $f \in   E$, one has $\int  \vert  f \vert  \, d\mu  =\int f \mathbf 1_{f>0}\, d\mu \,  -\, \int f \mathbf 1_{f<0}\, d\mu=\int f\, d\mu  '= \vert  \int f\, d\mu  ' \vert  $ where  $\mu  '=\mathbf 1_{f>0}\,\mu
-\mathbf 1_{f<0}\,\mu$ and $\mu  ' \in   { \mathcal K}-{ \mathcal K}$. Hence $ \vert   \vert   \vert  f \vert   \vert   \vert   \leq   \Vert f \Vert $.

The topology on   $E$ defined by the norm $ \Vert . \Vert $ is the topology of uniform convergence on ${ \mathcal K}-{ \mathcal K}$. Since ${ \mathcal K}-{ \mathcal K}$ is (convex, symmetric) vaguely compact (i.e. compact for $  \sigma  ({ \mathcal E},E)$) and generates ${ \mathcal E}$, Mackey's Theorem (see \cite{BOUEVT} or  \cite{SCHAE}) says that the topological dual of $(E, \Vert . \Vert )$ is ${ \mathcal E}$. For sake of completeness we recall the argument in our situation: by the bipolar theorem, the polar of the unit ball  $B_E$ in $(E, \Vert . \Vert )$ with respect to the duality  $(E;E^*)$ ($E^*$ is the algebraic dual) is ${ \mathcal K}-{ \mathcal K}$ (which -being $ \sigma  (E^*,E)$ compact- is closed for $ \sigma  (E^*,E)$  in $E^*$). One uses here that $B_E$ is by definition of the norm $ \Vert . \Vert $, the polar of ${ \mathcal K}-{ \mathcal K}$ for the duality $(E,E^*)$. But the polar of $B_E$ is also the unit ball of the dual of $(E, \Vert . \Vert )$. So the set of linear forms on $(E, \Vert . \Vert )$ whose norm is less than $1$  is exactly ${ \mathcal K}-{ \mathcal K}$ and the dual of  $(E, \Vert . \Vert )$ is ${ \mathcal E}$, equipped with the  gauge norm of ${ \mathcal K}-{ \mathcal K}$.
\medskip

In what follows $\nu $  denotes a nonnegative Borel measure in $X$ majorized by $\mathrm {Cap}$.

 \begin{proposition} \label{sci} The map $f \mapsto \int f^+\, d\nu  $ is lower semicontinuous in $(E, \Vert . \Vert )$.
\end{proposition}

We first notice the following  simple Markov-type inequality:
\begin{lemma} If $f \in   { \mathcal E}$ and $a>0$ then $${ \mathrm {Cap}}( \{  \vert  f \vert   \geq  a  \}) \leq  \vert   \vert   \vert  f \vert   \vert   \vert  / a .$$
\end{lemma}

{\sl Proof.} It suffices to note that  $  \mu   \{  \vert  f \vert   \geq a  \} \leq { \frac { 1  } {a}}(\int  \vert  f \vert  \, d\mu  ) \leq  { \frac { 1  } {a}}\vert   \vert   \vert  f \vert   \vert   \vert  $ for $\mu   \in   { \mathcal K}$. Taking the supremum over all $\mu   \in   { \mathcal K}$ the result follows. \qed

\bigskip

{\sl Proof of  proposition \ref {sci}.} Let $f_n \to  f$ in $E$. Since $ \vert  f_n^+-f^+ \vert   \leq  \vert  f_n-f \vert  $ we have $  \vert   \vert   \vert  f_n^+-f^+ \vert   \vert   \vert   \to  0$, and since $\vert  f_n^+\wedge f^+-f^+ \vert   \leq  \vert  f_n^+-f^+ \vert  $ we also have $ \vert   \vert   \vert  f_n^+\wedge f^+-f^+ \vert   \vert   \vert   \to  0$.

Now, given $ \varepsilon >0$ we have for every integer $n \geq 1$,
\begin{align}\int (f_n^+\wedge f^+)\, d\nu  &\geq \int f^+\, d\nu - \int_ {f_n^+\wedge f^+ \geq f^+- \varepsilon } (f^+-f_n^+\wedge f^+) \, d\nu  \nonumber
\\ &\hspace {6truemm} - \int_ {f_n^+\wedge f^+ < f^+- \varepsilon } (f^+-f_n^+\wedge f^+) \, d\nu \nonumber
\\ &\geq \int f^+\, d\nu - \varepsilon  \Vert \nu  \Vert - \Vert f \Vert _ \infty \,\,\nu (f_n^+ \leq f^+- \varepsilon ) \nonumber
\\ & \geq  \int f^+\, d\nu - \varepsilon  \Vert \nu  \Vert - \Vert f \Vert _ \infty \,\,{ \mathrm {Cap}} ( \{ f_n^+ \leq f^+- \varepsilon  \})\nonumber
\end{align}
using in the last line the inequality  $\nu  \leq {\mathrm {Cap}}$. By the preliminary remark  $${ \mathrm {Cap}} ( \{ f_n^+ \leq f^+- \varepsilon \} ) \leq  \vert   \vert   \vert  f_n^+\wedge f^+-f^+ \vert   \vert   \vert  / \varepsilon   \to  0 \mathrm{\ as\ } n \to   \infty .$$

Thus writing $ \int f_n^+\, d \nu  =\int (f_n^+\wedge f^+)\, d \nu  +\int (f_n^+-f_n^+\wedge f^+)\, d \nu  $ and observing  that $\int (f_n^+-f_n^+\wedge f^+)\, d \nu   \geq 0$ one sees that $\liminf \int f_n^+\, d \nu   \geq \int f^+\, d \nu   - \varepsilon  \Vert  \nu   \Vert $. Letting $ \varepsilon   \to    0$ the result follows: $\liminf \int f_n^+\, d \nu   \geq \int f^+\, d \nu $. \qed

We may now conclude the proof of Theorem \ref{rich}. Recall that $\nu  \in   {\mathfrak M}_+(X)$ and $\nu  \leq { \mathrm {Cap}}$.

 \begin{proposition} The measure $\nu $ is the upper envelope of the measures $\mu   \in   { \mathcal E}_+:= \cup _{n \geq 1} n\, { \mathcal K}$ such that $\mu   \leq \nu $. \end{proposition}

{\sl Proof.} The function  $ \Phi :f \to  \int f^+\, d \mu$ is convex, l.s.c.\ and homogeneous on $(E, \Vert . \Vert )$. The epigraph of $ \Phi$ is hence a closed convex cone in $E\times { \mathbb  R}$. By the Hahn-Banach separation theorem, $ \Phi $ is the upper envelope of the set of all continuous affine functions in $E$ that are less than $ \Phi $ in $E$. Now if $F=a+\ell$  is such a function with $\ell  \in   { \mathcal E}$ and  $a \in   { \mathbb  R}$, then   $\ell$ must be nondecreasing in $E$ (since $\ell$ is upper bounded in $E_-:=-E_+$ and hence $ \leq 0$ in $-E_+$). Since  for $f \in   E_+$, $a+t\ell(f)\leq tá\Phi (f)$ for all $t\geq0$, we have $\ell (f) \leq  \Phi (f)$  and
\begin{align} \int f\, d\nu &=\sup  \{ \ell(f)\, ;\,\ell  \in   E'    \; {\mathrm {such\ that\ }} \ \ell \leq \Phi    \; \}\nonumber
\\ &=\sup  \{\int f\,d\mu; \mu   \in   { \mathcal E}_+,\; \mu  \leq  \Phi  {\mathrm{\ in} \  } E_+  \}. \nonumber
\end{align}
 Thus $\nu$ is the smallest measure that majorizes $ \{ \mu \in    { \mathcal E}_+\,   \,;\, \mu   \leq \nu  \}$. This set being sup-stable in ${ \mathfrak M}_+(X)$,  $\nu  $ is the limit in $\mathfrak M_+(X)$ of an increasing sequence in ${ \mathcal E}_+ $. \qed

The following special cases is needed in the other parts of this paper. We use here the notations and general assumptions introduced in section 1. Denote for $a \geq 0$,
$${ \mathcal K}_{g,a}:= \{ \mu   \in   \mathfrak M_+(\partial \Omega )\,;\, \int _{ \Omega } g(x,\mathbb K_\mu  (x))\,  \psi (x)\, dx \leq a\, \}$$

\begin{corollary}\label{richExample} Let $\lambda   $ be a finite positive Borel measure in $\partial \Omega $ such that $\lambda  (K)=0$ whenever $K \subset   \partial \Omega $ is compact and ${ \mathcal K}_{g,a}$--null (i.e.\ $\mu  (K)=0$ for every $\mu   \in   { \mathcal K}_{g,a}$). Then there is an  increasing sequence $ \{\lambda   _n \}$  in $\mathfrak M_+(\partial \Omega )$ such that $\lambda   _n \in   n\, { \mathcal K}_{g,a}$ for $n \geq 1$  and $\lambda  =\lim_{n \to   \infty } \lambda _n $. \end{corollary}

{\sl Proof.} It suffices to note that ${ \mathcal K}_{g,a}$ is convex and weakly closed in ${ \mathcal M}_+(\partial \Omega) $. These properties follow from the convexity of $g(x,.)$, the continuity of $g$ and Fatou's theorem. $\square$

\begin{corollary}\label{richExample2} Suppose that $g$ satisfies the $ \Delta  _2$-condition and that $\mu   \in   { \mathfrak M}_+^{g, \Delta  }(\partial \Omega )$ vanishes on $\mathfrak M_+^{g, \Delta  }$--null sets. Then $\mu =\lim_{n \to   \infty } \mu_n $ for some increasing sequence $ \{ \mu  _n \}$ in $\mathfrak M_+^{g, \Delta}(\partial \Omega )$.
\end{corollary}
{\sl Proof.}  It suffices to apply the previous result with $a=1$ and observe that by the assumptions the cone $\mathfrak M_+^{g, \Delta  }(\partial  \Omega )$ is generated by ${ \mathcal K}_{g,1}$. $\square$

\section{Some known results on Schr\"odinger equations}

 In this section, we present some  notions and results related to the Schr\"odinger equation (\ref{Sch}) that are needed to obtain our results. Here $V$ is a nonnegative continuous function in $\Omega  $ satisfying condition  \eqref{Vcond'}, namely $V \in   { \mathcal V}(\Omega ,\overline  a)$ with $\bar a \geq 0$.

 We start with some definitions.

 {\sl $L^V$ harmonicity, $L^V$ superharmonicity.} By an {\sl $L^V$ harmonic function} in $\Omega $ we mean a continuous function in $\Omega $ that solves  $L^V(u)=0$ in the sense of distributions. Then $u \in   W_{\rm loc}^{2,p}(\Omega )$ for every $p< \infty $ and $L^V(u)=0$ a.e. As well known  every $u \in   L^{1}_{\rm loc}(\Omega )$ solving $L^V(u)=0$  in the distribution sense in $\Omega $ can be viewed as an $L^V$ harmonic function.

 It is also known that if $u \in   L^1_{\rm loc} (\Omega )$ is a supersolution of (\ref{Sch}) (i.e.\ $L^V(u)  \geq 0$ in the weak sense) then the limit  $ \tilde u(x):=\lim_{r \to  0} \int _{B(x,r)}\, u(y)\, dy/ \vert  B(x,r) \vert   \in    (- \infty ,+ \infty ]$ exists at every $x \in   \Omega $, $ \tilde u$ is l.s.c.\  and $L^V$-superharmonic w.r.\ to the harmonic sheaf of $L^V$-harmonic functions (in the sense of \cite{brelot}). Conversely every $L^V$-superharmonic function in a domain $\omega  \subset   \Omega $ which is not the constant $+ \infty $ is a locally integrable supersolution in $\omega $.  In the sequel we will only need to consider continuous $L^V$ superharmonic functions.

  {\sl $L^V$ potentials.} A non-negative, $L^V$ superharmonic function   in $\Omega $ is an $L^V$ potential if it does not have a  positive $L^V$ harmonic minorant in $\Omega $.

  {\sl Integral representations of nonnegative harmonic functions.} By the results in \cite{Annals87} (see also the presentation in \cite{An-MZ}) there exists a continuous kernel $\mathbb K^V: \partial \Omega \times \Omega \ni ( \zeta  ,x)  \mapsto  \mathbb K_ \zeta  (x)  \in   (0,+ \infty )$ such that  $\mathbb K_ \zeta (x) =\lim_{z \to   \zeta  } \mathbb  G(z  ,x)/\mathbb  G(z  ,x_0)$ for $ \zeta   \in   \partial \Omega $, $x \in   \Omega $ (thus $\mathbb K_ \zeta  (x_0)=1$). Moreover,
 every positive $L^V$-- harmonic function $u$ in $\Omega $ can be uniquely written in the form  $u=\mathbb K^V_\mu  :=\int \mathbb K_ \zeta  (.)\, d\mu  ( \zeta  )$ where $\mu   \in   \mathfrak M_+(\partial \Omega )$. The function $\mathbb K$ is the $L^V$--Martin kernel for $\Omega $ with normalization at $x_0$ and if as above  $u=\mathbb K_\mu  ^V$,   $\mu  $ will be called  the $L^V$ boundary measure of $u$.

  {\sl $L^V$ moderate solutions and supersolutions.} The next definition and statements are parallel to Definition \ref{gmod-sol}, Proposition \ref{gmodprop} and Remark \ref{rksuper}.

 \begin{definition}\label{LVnotions}
 An $L^V$ harmonic function $u$ is \emph{$L^V$-moderate} (or $V$-moderate) if $uV\in L^1_\psi(\Omega )$ that is: $\int _{\Omega }  \,\vert  u \vert  \,V\,  \psi \, dx < \infty $.
 Similarly a positive $L^V$ superharmonic function $u$ is $L^V$-moderate if $uV\in L^1_\psi(\Omega )$.
 \end{definition}

\begin{proposition} Let $u$ be nonnegative and $L^V$ harmonic in $\Omega $. Then the following are equivalent: (i) $u$ is $L^V$-moderate, (ii) $u$ has a $m$-boundary trace in $\partial \Omega $, (iii) $u$ admits a harmonic majorant in $\Omega $. \end{proposition}

\proof  (i) implies (ii): see proof of Lemma \ref{A0}.

(ii) implies (iii): the harmonic function with the same m-boundary trace dominates $u$.

(iii) implies (i):   The infimum of all superharmonic majorants, say $\bar u$,  is superharmonic. Therefore $v:=\bar u-u$ is a non-negative superharmonic function. Thus $v$ has an m-boundary trace $\tau\geq 0$.  Since $v':=\bar u-\mathbb K_\tau\geq u$ and $v'$ is superharmonic it follows that $v=v'$ so that $\tau=0$. It follows that $u$ has an m-boundary trace equal to $\mathrm{tr}\,\bar u$.
\qed

 The following lemma is proved in \cite{MM-mod} in the case of $C^2$ domains.

 \begin{lemma} \label{A0} Let $w$ be a positive  $L^V$ superharmonic function in $\Omega $. If $w$ is $L^V $ moderate
then  it possesses an  $m$-boundary trace. The largest $L^V$ harmonic function dominated by $w$ is    $L^V $  moderate and its  $m$-boundary trace equals $ {\mathrm tr}\,w$.
\end{lemma}

\proof The function $ \widehat w:=w+ { \mathbb G} [wV]$ is  positive superharmonic in $\Omega $.  Consequently it  has an  $m$-boundary trace $\nu\in {\mathfrak M}( \partial \Omega )$ and since $ {\mathrm tr}\, { \mathbb G} [wV]=0$ it follows that $\nu= {\mathrm tr}\,w$.

Also, because $ \widehat w:=w+ { \mathbb G} [wV]$ is  positive superharmonic, $-  \psi  \Delta   \widehat w \in   {{\mathfrak M}}_+(\Omega )$, i.e.\ $ -\psi ( \Delta  w-Vw)=- \psi  L^V (w)\in   {{\mathfrak M}}_+(\Omega )$.
Since $0 \leq  { \mathbb G} ^V(- \Delta  w+Vw) \leq  { \mathbb G} (- \Delta  w+Vw)\not \equiv + \infty $ this implies that $v:= { \mathbb G} ^V(- \Delta  w+Vw)$ has a vanishing boundary trace. Finally $w-v$ is the largest $L^V$ harmonic function dominated by $w$ and $\mathrm{tr}(w-v)= \mathrm{tr}\,w$.
 \qed

{\sl A Fatou-Doob-Na\"\i m type result.} For easy  reference we state some additional potential theoretic results   that are used in the sequel and which extend well-known results of Fatou, Doob and Na\" \i m. We refer the reader to \cite{Annals87}. See also  \cite{An-SLN} or \cite{An-MZ}. Recall that the \emph{Riesz decomposition theorem} for $L^V$ says that a non-negative $L^V$ superharmonic function $u$ in $\Omega $, admits a greatest  $L^V$-harmonic function $h$ dominated by $u$ and that $u-h$ is a potential.

 \begin{theorem}\label{Fatou}  a) If $u,v$ are positive $L^V$ harmonic functions in $\Omega  $ then,
$$\lim_{x\to \zeta  } \frac{u(x)}{v(x)}=\frac{d\mu_u}{d\mu_v}( \zeta  )   \quad   L^V-finely, \;\mu_v - \text{a.e. on $ \partial \Omega $}$$
where $\mu_u$ and $\mu_v$ are the $L^V$ boundary measures of $u$ and $v$ respectively and the term on the  right hand side  denotes the Radon-Nikodym derivative.\footnote{See the comments on $L^V$-fine convergence, following the theorem.}

(b) If $p$ is an $L^V$ potential and $v$ is positive $L^V$ harmonic function:
$$\lim_{x\to \zeta  } p(x)/v(x)=0  \quad   L^V-finely, \;\mu_v - \text{a.e. on $\partial \Omega $}$$
where $\mu_v$ is the $L^V$ boundary measure of $v$.
\end{theorem}

For the definition of $L^V$ fine convergence see e.g.\ \cite{An-SLN} (where proofs of the above statements are provided).  It is useful and handy to note that if $f$ is a non-negative function in $\Omega $ satisfying the strong {\sl Harnack property} in $\Omega $ then (see \cite{An-MZ} p. 409).
\begin{equation}\label{fine-nt}
\lim_{x\to \zeta  } f(x)=a  \quad    L^V-finely \; \Rightarrow  \;  f(x)  \to  a\  as \  x \to   \zeta  \ non-tangentially.
\end{equation}
Here we say that a nonnegative function $f$ in $\Omega $ satisfies the strong (uniform) Harnack property in $\Omega $ if for every $ \varepsilon >0$ there is a positive $  \eta $ such that
\begin{equation}\label{HP}
 \vert  f(x)-f(y) \vert   \leq  \varepsilon f(x) {\rm \  whenever\ } x,\,y \in   \Omega  {\rm \ and \ }  \vert  x-y \vert   \leq   \eta  d(x;\partial \Omega )\end{equation}
Nonnegative solutions $u$ of (\ref{Sch}) or (\ref{geq})
 have this property. If $f$ and $g$ satisfy the Harnack property and $g$ is $>0$ in $\Omega $ then the ratio $f/g$ also satisfies the strong Harnack property in $\Omega $. In the present paper  Theorem \ref{Fatou}  will be applied only in cases where $p,u,v$ satisfy the Harnack property in $\Omega $. Therefore in our applications we may replace ``$L^V$ fine convergence''  by ``n.t. convergence'' (where n.t. is an abbreviation for non-tangential).

{\sl $L^V$--regularity of boundary points.} For the next definition note that for $ \zeta  \in \partial \Omega $ the harmonic function $\mathbb K _ \zeta  $ is $L^V$ superharmonic. Therefore, by the Riesz decomposition theorem, $\mathbb K _\zeta$ can be uniquely represented in the form
\begin{equation}\label{tlKV}
 \mathbb K _\zeta=\tilde { \mathbb K }_\zeta^V+p_\zeta
\end{equation}
where $\tilde  { \mathbb K }_\zeta^V$ is the largest $L^V$ harmonic function dominated by $  \mathbb K _\zeta$ and $p_\zeta$ is an $L^V$-potential. ( $\tilde  { \mathbb K }_\zeta^V$ may be the zero function.) $ { \mathbb K }^V_\zeta$ is the unique $L^V$ harmonic function that vanishes on $\partial\Omega\setminus\{\zeta\}$ and is normalized by $ { \mathbb K }^V_\zeta(x_0)=1$. Therefore,
\begin{equation}\label{c_V}
\tilde  { \mathbb K }_\zeta^V:=c_V(\zeta)  { \mathbb K }_\zeta^V,
\end{equation}
for some $c_V( \zeta  )\geq 0$. Note that $c_V( \zeta  ) \in [0,1] $.

\medskip

 \begin{definition}\label{V-regular}  (see  \cite{An-App}, \cite{An-MZ})  A point $ \zeta  \in \partial \Omega $ is $L^V$ \emph{regular} if $ { \mathbb K }_ \zeta  $ is not an $L^V$ potential, i.e.\ if $c_V( \zeta  )>0$. The point $ \zeta  $ is $L^V$ singular if it is not regular i.e.\ if $c_V( \zeta  )=0$.

  The set of $L^V$ regular points is denoted by $ { \mathcal R}^V$ while the set of $L^V$ singular points is denoted by $ { \mathcal S}^V$.
\end{definition}

The next facts are taken from  \cite{An-App} (see also \cite{An-MZ}).  Let $ \zeta   \in   \partial \Omega $. Let $\nu $ be a unit vector  in $   { \mathbb R} ^N$ be such that $\limsup_{x \in   \partial \Omega ,x \to   \zeta  }  {\frac  {\langle x- \zeta  ,\nu  \rangle } { \vert  x- \zeta   \vert  }}   <1$, i.e.\, $\nu $ is a pseudo-inner normal at $ \zeta  $.
 \begin{proposition}\label{sing-reg} (i)  $ \zeta   \in    { \mathcal S}^V$ if and only if
$  { \mathbb K }_ \zeta  (x)/ { \mathbb K }^V_ \zeta  (x)\to 0  \quad   \text{as $x\to \zeta  $ n.t.}$,

(ii) $ \zeta   \in   { \mathcal R}^V$ if and only if $\liminf_{t \to  0}  { \mathbb K }_ \zeta  (x_0+t\nu )/ { \mathbb K }^V_ \zeta  (x_0+t\nu)>0$,

(iii)  $ \zeta   \in    { \mathcal S}^V$ if and only if
$ {  \mathbb G } ^V_{x_0}  (x)/ {  \mathbb G } _ {x_0}  (x)\to 0$  as $x\to \zeta$ n.t.,

(iv) $ \zeta    \in   { \mathcal R}^V$ if and only if $\liminf_{t   \to     0  } {\frac  { {  \mathbb G } ^V_{x_0}  ( \zeta  +t\nu ) } {{  \mathbb G } _ {x_0}  ( \zeta  +t\nu}}  )>0$.

\end{proposition}

\medskip

The next useful fact connecting the $L^V$ boundary measure and the $m$-boundary trace for moderate $ \geq 0$ solutions  is proved in \cite{MM-mod} in the case of $C^2$ domains.

 \begin{proposition} \label{MBLB} If $u$ is a positive    $L^V $  moderate harmonic function then its m-boundary trace and its    $L^V $  boundary measure are mutually absolutely continuous. \end{proposition}

\noindent\textit{Remark.}\hskip 2mm Using Lemma \ref{A0} this is proved in the same way as Proposition 3.8 of \cite{MM-mod}. An alternative argument is given below.

\proof  Suppose $ {   \mathbb K } _\mu  ^V$ is moderate ($\mu   \in   {\mathfrak M}_+(\partial \Omega )$) and let $\lambda $ be its $m$-boundary trace on $\partial \Omega $. If $ \zeta  $ is regular then, by  \eqref{tlKV} and  \eqref{c_V},
 $$c_V( \zeta  )^{-1} {  { \mathbb K }} _  \zeta  =  {  { \mathbb K }} _ \zeta  ^V+p'_  \zeta  $$  and  $p'_ \zeta  = {  \mathbb G }  (V {   \mathbb K } _ \zeta  ^V)$. It follows that $ {   \mathbb K } _{c^{-1}_V\mu}  = {   \mathbb K } ^V_{\mu   } + {\rm potential}$ and $c_V^{-1}\mu  =\lambda $.  \qed

The following somewhat related fact will be useful later.

 \begin{lemma} \label{KVK} Suppose $V \in   { \mathcal V}(\bar a, \Omega )$, $\mu  $, $\nu \in   \mathfrak M_+(\partial \Omega )$ and assume that $ \mathbb K _{\nu}^V \leq  \mathbb K  _\mu$ in $\Omega $. Then for every Borel set $A \subset   \partial \Omega $, $$ \mathbb K _{\mathbf 1_A\, \nu}^V \leq  \mathbb K _{\mathbf 1_A\, \mu}$$

 $\mathbb K^V_\nu$ is $L^V$ moderate (see Proposition 3.2) and therefore possesses a boundary trace. Moreover, for every Borel set $A \subset   \partial \Omega $,
   $${\mathrm tr}( \mathbb K _{\mathbf 1_A\,\nu} ^V)=\mathbf 1_A\, \mathrm{tr}\,(\mathbb K _{\nu}^V).$$
   \end{lemma}

       {\sl Proof.}
For the first claim, it suffices -by a standard approximation argument using the regularity of finite Borel measures- to show that $ \mathbb K _{\mathbf 1_F\, \nu}^V \leq  \mathbb K _{\mathbf 1_U\, \mu}$ whenever $ F$ is compact, $U$ is open in $\partial \Omega $ and $F\subset   U$. To prove this observe that $\inf \{ \mathbb K _{\mathbf 1_U\, \mu}-  \mathbb K _{\mathbf 1_ F\, \nu}^V ,0 \}$ is a classical superharmonic function in $\Omega $ which has a nonnegative inferior limit at every point of $\partial \Omega $ (use the fact that $\mathbb K _ \lambda ^V$ goes to zero at every $ \zeta   \in   \partial \Omega \setminus \text{supp} \lambda $). Thus it is nonnegative and the claim follows.

       Put $\lambda= \mathrm{tr}\,\mathbb K _{\nu}^V$.
       Obviously $ \mathbb K _{\nu}^V \leq  \mathbb K  _\lambda$ in $\Omega $. Therefore, if $ A$ is a Borel subset of $\partial \Omega $, $ \mathbb K _{\mathbf 1_A\, \lambda}- \mathbb K ^V_{\mathbf 1_ A\, \nu}$ is a non-negative (classical) superharmonic function.  Again, by the first part of the lemma (applied to $A'=\partial \Omega \setminus  A$) we have
       $$ \mathbb K _{\mathbf 1_ A\, \lambda}- \mathbb K ^V_{\mathbf 1_ A\, \nu}\leq  \mathbb K _{\lambda}- \mathbb K ^V_{ \nu}.$$
       Since  $  \mathbb K  _\lambda- \mathbb K _{\nu}^V$ is a potential in $\Omega $, the nonnegative superharmonic function  $ \mathbb K _{\mathbf 1_ A\, \lambda}- \mathbb K ^V_{\mathbf 1_ A\, \nu}$ is also a potential and its $m$ boundary trace is zero. \qed

\bigskip

{\sl Two regularity criteria.} The next result was proved in \cite{An-MZ} and \cite{MM-mod}. (In \cite[Proposition 3.9]{MM-mod}  the  result was proved for $C^2$ domains. Using Lemma \ref{A0}, the  proof applies equally well to Lipschitz domains.)

\begin{proposition}\label{A2}  For every $ \zeta  \in  \partial \Omega $
\begin{equation}\label{MM3.1}
\int_\Omega     {\mathbb K }_ \zeta  ^V V\psi\,dx<\infty \iff  \zeta   \in  { \mathcal R}^V.
\end{equation}
Furthermore, if $u$ is a positive $L^V$ moderate harmonic function and  $\nu= {\mathrm tr} _{ \partial \Omega }u$ then,
\begin{equation}\label{MM3.1ii}
 \nu( { \mathcal S}^V)=0.
\end{equation}
\end{proposition}

 \medskip

Finally we mention the following characterization of    $L^V $  regularity which is the main result in \cite{An-MZ}.
 \begin{proposition}\label{A3} Let $ \zeta  \in  \partial \Omega $. Then
\begin{equation}\label{Ancona1}
 \zeta  \in  { \mathcal R}^V\iff   \int_\Omega    { \mathbb K }_ \zeta  \,V\psi\,dx<\infty.
\end{equation}
\end{proposition}

       \section{The Boundary Harnack Principle and related auxiliary results}\label{s:BHP}

Let   $r,\,  \rho $ be positive numbers and let $f:{ \mathbb  R}^ { N-1} \to  { \mathbb  R}$ be Lipschitz with a Lipschitz constant $ \leq { \frac {   \rho } {10\, r }}$ and such that $f(0)=0$.  We set
       $$\omega _f(r, \rho ):= \{  \xi  =( \xi  ', \xi  _N) \in   { \mathbb  R}^{N-1}\times { \mathbb  R}\, ;\, \vert   \xi  ' \vert   < r,\, f( \xi  ')< \xi  _N< \rho  \}.$$
  This region will be called  the   standard Lipschitz domain of height $ \rho $, radius $r$ and defining function $f$. We denote $A= A(\rho ):=(0,\dots, 0,{ \frac {   \rho } {2 }})$,  $A'= A'( \rho ):={ \frac {  3} {4 }}\,A$.

      We assume at first that $\Omega $ is a region in ${ \mathbb  R}^N$ whose intersection with the cylinder $ T(r, \rho ):= \{ \xi  =( \xi  ', \xi  _N) \in   { \mathbb  R}^{N-1}\times { \mathbb  R}\, ;\, \vert   \xi  ' \vert   < r,\, - \rho < \xi  _N< \rho  \}$ is  $\omega =\omega _f(r, \rho )$.

     We fix some $V \in   { \mathcal V}(\bar a, \Omega)$, $\bar a>0$. As before,  $\mathbb G ^V$ denote the Green's function for $L^V$ in $\Omega $ and $\mathbb G =\mathbb G ^0$. The following boundary Harnack principle for our operators $L^V$ is proved in \cite{Annals87} (see also a presentation in \cite{An-MZ}).

       \begin{theorem} \label{HBI}
      There is a constant $c$ depending only on $N,\, \bar a$ and ${ \frac {   \rho } { r}}$ such that whenever $u$ is a positive $L^V$ harmonic function  in $\omega$  that vanishes continuously in $\partial \omega \cap T(r, \rho )$
 then
\begin{equation}\label{u/v+}
\hspace{2truemm}  c^{-1}r^{N-2}\, \mathbb G ^{V}(x,A')\leq \frac{u(x)}{u(A)}\leq c\, r^{N-2} \,\mathbb G ^{V}(x,A'), \quad  \hspace{7truemm} \forall x\in \Omega  \cap \overline T({ \frac {  r} {2 }}; { \frac {   \rho } {2 }})
\end{equation}
In particular,  for any pair $u,v $ of  positive $L^V$ harmonic functions in  $ \omega $ that vanish on $\partial \omega  \cap T(r, \rho )$:
\begin{equation}\label{u/v}
 u(x)/v(x)\le Cu(A)/v(A), \quad \quad  \forall x\in \omega _f(r/2, \rho /2))
\end{equation}
where $C=c^2$.
 \end{theorem}

{\sl Remark.} Suppose that  $u, \, v$ is a pair of nonnegative $L^V$ harmonic functions  in $\omega '=\omega \setminus \overline T(r/4, \rho /4)$ that  both vanish on $\partial \omega  \cap (T(r, \rho )\setminus \overline T(r/4, \rho /4))$. Using  Theorem \ref{HBI}, it is easy to see that, for  $x$  restricted to $\omega   \cap \partial T(r/2, \rho /2)$,  (\ref{u/v}) holds again -with another constant $C=C(\overline a,r/ \rho ,N)$-.

Actually we need a result (Proposition \ref{mainSec4} below)  similar  to  \eqref{u/v} in Theorem \ref{HBI} above but  for the ratio of solutions of two \emph{different equations},  namely the equations $L^V(u)=0$ and $ \Delta  v=0$. Various  results of this type are known, in particular it can be obtained  as a consequence of an assumption on the growth of $V$ necessarily stronger than \eqref{Vcond'}  (see e.g. \cite{An-CR82}, \cite{An-Green97} and the references therein). Here  we will derive estimates of  the ratio of $\mathbb K _\zeta/\mathbb K ^V_\zeta$ in $\Omega$ uniform w.r. to $\zeta$ in certain subsets of ${ \mathcal R}^V$   (see Proposition \ref{mainSec4} below).
To this end we use the boundary Harnack principle  Theorem \ref{HBI} together with its following well-known consequence (see e.g.\   \cite[Lemma 3.5]{An-MZ}).
It will be assumed in the next four statements that the reference point $x_0$ is in $\Omega \setminus \overline  T(r, \rho )$.

 \begin{proposition}\label{FMG}
 There exists a constant $C$  such that for all  $x= (0,\dots, t)$,  $\vert  t \vert   \leq { \frac {  3} {4 }} \rho $, \begin{equation}\label{An-3.5}
  C^{-1}t^{2-N}\leq  { \mathbb K }^V_0  (x)\mathbb G ^V_{x_0}(x)\leq C\,t^{2-N}.
\end{equation}
and $C$ can be chosen depending only on $\bar a$, ${ \frac {   \rho } {r }}$ and $N$.
\end{proposition}

We now formulate conditions  in terms of the Green kernels $\mathbb G _{x_0}$ and $\mathbb G ^V_{x_0}$ for the Martin kernels $ { \mathbb K }_ \zeta  $ and $ { \mathbb K }^V_ \zeta  $ to be equivalent.

            \begin{proposition}
 Suppose that for some constant $C  \geq 1$ it holds that
$$\mathbb G (x_0,y) \leq C\, \mathbb G ^V(x_0,y)\ \   for\ all\  y \in   \omega _f(r, \rho ).$$
Then there is a constant $C_1\geq 1$ depending on $C$, ${ \frac {   \rho } {r }} $, $N$ and $\bar a$ such that for all $ \zeta   \in  \partial  \Omega  \cap \overline T(r/8, \rho/8 )$ and $x \in \Omega  \cap \overline T(r/4, \rho/4 )$
$$C_1^{-1}  { \mathbb K }_ \zeta  (x) \leq  {\mathbb K }_ \zeta  ^V(x) \leq C_1\,  { \mathbb K }_ \zeta  (x) $$
\end{proposition}
Recall that $\mathbb G^V \leq \mathbb G$ always holds since $V \geq 0$.

\noindent{\sl Proof.}  It suffices to consider the case $ \zeta  =0$. Choose $0<t \leq 1/4$ such that $x \in   \partial( T(tr, t \rho ))$. It follows from the remark following Theorem \ref{HBI} that
\begin{align} \label{GVBH}       C_1^{-1} \,\mathbb G ^V(z,A'(t \rho )))\,\mathbb G ^V(A'(t \rho ),x) \,(tr)^{N-2} & \leq \mathbb G ^V(z,x)
   \\  \nopagebreak \nonumber  \leq  C_1\,\mathbb G ^V&(z,A'(t \rho )))\,\mathbb G ^V(A'(t \rho ),x) \,(tr)^{N-2}
   \end{align}
when $z \in  \Omega  \cap  \partial T(2tr,2t \rho )$. By the maximum principle for $L^V$ these relations extend  to all $z \in  \Omega  \setminus T(2tr,2t \rho )$ and thus hold at $z=x_0$. Using  \eqref{GVBH} at $z=x_0$, together with their analogues for $V=0$ and the hypothesis, one obtains  that $\mathbb G ^V(A'(t \rho ),x)$ and $\mathbb G (A'(t \rho ),x)$ are uniformly comparable. By   Theorem \ref{HBI}, $ { \mathbb K }_ \zeta  ^V(x)$ is uniformly comparable to $ { \mathbb K }_ \zeta  ^V(A'(t \rho ))\, \mathbb G ^V(A'(t \rho );x)(tr)^{N-2}$, hence also to the ratio $\, \mathbb G ^V(A'(t \rho );x)/\mathbb G ^V(x_0,A'(t \rho ))$ (using Proposition \ref{FMG}). The proposition follows. $\square$

 \begin{corollary} \label{cor4}Assume that $\mathbb G (x_0,y) \leq C\, \mathbb G ^V(x_0,y)$ for  $y \in   \Omega $ (not only in $\Omega \cap T(r, \rho )$ as above). Then there is a constant $C_2=C_2(\Omega ,r, \rho,$ $x_0,\bar a) \geq 1$ such that for  $ \zeta   \in  \partial \Omega  \cap  T(r/8, \rho /8)$ and  $x \in   \Omega $:
 \begin{align}   \nonumber    C_2^{-1}  { \mathbb K }_ \zeta  (x)  \leq  { \mathbb K}_ \zeta   ^V(x)
 \leq  C_2\, { \mathbb K}_ \zeta   (x)
   \end{align}

 \end{corollary}
 {\sl Proof.} In fact it follows from Theorem \ref{HBI} that $ {\mathbb K}_ \zeta^V  (x)$ is -for $x \in   \Omega \setminus T({ \frac {  r} {2 }};{ \frac {   \rho } {2 }})$ and given $\Omega $, $x_0$, $r,\, \rho $, $\bar a$- uniformly equivalent to $\mathbb G ^V(A'({ \frac {   \rho } {2 }}),x)$ (it suffices to consider the case $x \in   \partial T({ \frac {  r} {2 }};{ \frac {   \rho } {2 }})$).   Using Harnack inequalities (for $\Omega \ni y \mapsto \mathbb G ^V(y,x)$) one then sees that $ { \mathbb K }_ \zeta^V  (x)$ is uniformly equivalent to  $1\wedge \mathbb G ^V(x_0, x)$ --under the same assumptions on $x$ and $ \zeta  $--. Using this also for $V=0$ together with the hypothesis the result follows for $x \in   \Omega \setminus T({ \frac {  r} {2 }};{ \frac {   \rho } {2 }})$. $\square$

 \begin{lemma} \label{lemcomp} Suppose that  $0 \leq V \leq { \frac {  C} {r^2 }}$ in $\omega =\omega _f(r, \rho )$. Then, with $A=A({ \frac {   \rho } {2 }})$ and for some constant $C_1=C_1({ \frac {  r} { \rho  }}, N,C)$
$$\mathbb G _{x_0}^V  \leq  \mathbb G _{x_0} \leq C\,{\frac  { \mathbb G (x_0,A)} {\mathbb G ^V(x_0,A)}} \, \mathbb G _{x_0}^V \  \mathrm{in} \ \omega_f({ \frac {  9r} {10 }};{ \frac {  9 \rho } {10 }})$$
\end{lemma}

 We remark here  that by homogeneity, one may assume that $r=1$ and $ \rho $ is fixed. Some general similar comparability results (w.r.\ to two different operators) appear in \cite{An-Green97}.

 {\sl Proof.} Using the Harnack boundary principle Theorem \ref{HBI}, it suffices to show that the Green's functions $g$, $g^V$ for $L^0$ and $L^V$ in $\omega $ are such that $g(A,.) \leq c\, g^V(A,.) $ in $\omega $ with a constant $c=c(N,{\frac  { r} { \rho }} , C)$. A proof of the later fact is provided in \cite{An-CR82}. In fact it would suffice for the proof to use \cite[Proposition 3]{An-CR82}. \qed

From now on $\Omega $ denotes a general bounded Lipschitz domain in ${ \mathbb R} ^N$.  By definition for each $P \in   \partial \Omega $ there is an isometry $I:{ \mathbb  R}^N \to  { \mathbb  R}^N$ and a cylinder $T(r, \rho )$ such that $I(P)=0$ and $I(\Omega ) \cap T(r, \rho )$ is a standard Lipschitz domain  in ${ \mathbb  R}^N$ with height $ \rho $ and radius $ r$. So the previous propositions can be applied to $\Omega _0=I(\Omega )$ and $T(r, \rho )$.

 If $ \zeta   \in   \partial \Omega $, $x \in   \Omega $ and $0< \theta  <{ \frac {   \pi  } {2 }}$, we will say that $[\zeta  ,x]:= \{  \zeta  +tx\,;\, 0<t \leq 1\,  \}$ is inner $ \theta  $-non-tangential at $ \zeta  $ (w.r.\  to $\Omega $) if the truncated cone $C( \zeta  ,x,  \theta   ):= \{ z\,;\,    \cos  \theta  \; \vert z- \zeta   \vert \,  \vert  x- \zeta   \vert  < \langle  z- \zeta  ,x- \zeta   \rangle\,< \vert  x- \zeta   \vert  ^2    \} $ is contained in $\Omega $.

       In the next statement  $ F$  denotes a
compact subset of $\partial \Omega $    such that $V(x) \leq { \frac {  \bar a} {d(x, F)^2 }}
$ for $x \in   \Omega $. As before a reference point $x_0 \in   \Omega $ is fixed.

       \begin{proposition}\label{mainSec4} Let $0< \theta  < { \frac {   \pi  } {2 }}$, $c_1,\, c_2>0$ and assume that   for every $ \zeta   \in  F$ there is a inner $ \theta  $-non-tangential  segment  $[ \zeta  ,b_ \zeta  ]$ at $ \zeta  $ (w.r.\ to $\Omega $ )   with $ \vert  b_ \zeta  - \zeta   \vert   \geq c_2$ and
       $$\hspace{8truemm}  { \mathbb K }^V_ \zeta (x) \leq c_1  { \mathbb K }_ \zeta (x), \hspace{10mm}   \forall    \zeta   \in    F,\;\; \forall   x \in   ( \zeta  ,b_ \zeta  ]$$ (thus $F \subset   { \mathcal R}^V$).
Then there exists a  constant $\bar c \geq 1$ depending only on $\bar a$, $ \theta  $, $c_1$, $c_2$ and  $\Omega $ such that, for all $ \zeta   \in   \partial \Omega $ and all $x \in   \Omega $
$$\hspace{15mm}
 \bar c^{\,-1}\leq \frac{ { \mathbb K }_ \zeta  (x)}{ { \mathbb K }_ \zeta  ^V(x)}\leq \bar c \ \ \  $$
       \end{proposition}

       {\sl Proof.} By  Proposition \ref{FMG} the assumption means that for some constant $c \geq 1$, $\mathbb G (x_0, y) \leq c\, \mathbb G ^V(x_0,y)$ for all $y \in  \bigcup    ( \zeta  ,b_ \zeta  ]$. By Harnack inequalities and the definition of a Lipschitz domain this means that for every given $\kappa >0$, it holds that for some constant $c'=c'(c_1,c_2, \bar a, \Omega )>0$,
       $$\mathbb G (x_0,y) \leq c'\, \mathbb G ^V(x_0,y)$$
        for $y \in  U(\Omega , \kappa):=\Omega  \cap  \{ x\,;\, d(x, F) \leq \kappa d(x,\partial \Omega )\, \}$. Using this  for $ \kappa >0$ small enough,  Lemma \ref{lemcomp}  shows that $\mathbb G (x_0,y) \leq c\, \mathbb G ^V(x_0,y)$   whenever $y \in  \Omega \setminus U(\Omega ,2\kappa)$.  The proposition follows using Corollary \ref{cor4}  above. $\square$

         We close this section with another auxiliary result. Here we assume that the region
$\Omega $ is Lipschitz and bounded,  that $g:\Omega \times { \mathbb  R}_+ \to  { \mathbb  R}_+$ is continuous and that $g(x,t)/t$ is non decreasing in $t$, $t>0$.

\begin{proposition} \label{aux2}Let $u$ be a nonnegative solution of the equation
$$ \Delta  u(x)=g(x,u(x))$$ and let $F \subset   \partial \Omega $ be closed. Assume further that $g(x,u(x)) \leq { \frac {  c_0} {( \delta (x))^2 }}\, u(x)$ in $\Omega $ for some constant $c_0 \geq 1$. Then there is a largest solution $v$ dominated by $u$
 in $\Omega $ and vanishing on $\partial \Omega \setminus F$. This solution is also the largest nonnegative (say continuous) subsolution dominated by $u$ that vanishes on $\partial \Omega \setminus F$. \end{proposition}

{\sl Proof.} a) If $v$ is a nonnegative solution dominated by $u$ in $\Omega $, then $ \Delta  v=V\, v$ for some potential $V=V_v \in   { \mathcal V}(\Omega ,c_0)$ and $v= \mathbb K ^{V}_\mu$ for some $\mu \in   \mathfrak M_+(\partial \Omega )$.

Suppose that $v=0$ on $T(r, \rho ) \cap  \partial \Omega $ for some cylinder $T(r, \rho )$ such that $\omega (r, \rho )=\Omega  \cap T(r, \rho )$ is a standard Lipschitz domain of height $ \rho $ and radius $r$. Then by the boundary Harnack principle applied to the $L^V$-harmonic function $v$,
$$v(x) \leq C\, v(A') r^{N-2} \mathbb G _{A'}^{V}(x) \leq C u(A') r^{N-2} \mathbb G _{A'}(x), \quad x \in   \omega (r/2, \rho /2)$$
with $C$ depending only on ${ \frac {   \rho } {r }}$, $N$ and $c_0$. Here $A'=(0,{ \frac {  3} {4 }} \rho )$.

This shows that the family of all non-negative solutions of $ \Delta  u-g(x,u)=0$ dominated by $u$ and vanishing on $\partial \Omega \setminus F$  must uniformly  vanish at each point $ \zeta   \in   \partial \Omega \setminus F$

 {\sl Remark.} Recall that since $s=1$ is $L^{V_v}$ superharmonic, the Fatou-Doob  property available here implies that $\mu  (\partial \omega  \cap T(r, \rho ))=0$ (see Remark 2.7 in \cite{An-MZ}).

b) Let again $T(r, \rho )$ be a cylinder such that  $\omega (r, \rho )=\Omega  \cap T(r, \rho )$ is a standard Lipschitz domain of height $ \rho $ and radius $r$ and let now $v$ be finite non-negative lsc  function in $\overline \omega (r, \rho )$ which is $ \leq u$ in $\Omega  \cap  \overline \omega (r, \rho ) $, vanishes continuously on $\partial \Omega  \cap \overline T(r, \rho )$ and which is a subsolution in $\omega $. Since e.g.\ $ \Delta  v \geq 0$ we may assume after a modification on a negligible set that $v$ is subharmonic.

Setting again $V=g(x,v(x))/v(x)$ ($=0$ if $v(x)=0$), $v$ is clearly $L^{V_v}$ subharmonic in $\omega $ and $V_v \in   { \mathcal V}(\omega ,c_0)$.

Fix $ \varepsilon >0$ small. The function $(v- \varepsilon )_+$ is again a subsolution on $\omega _{ \varepsilon '}:=T(3r/4, \rho ) \cap    (\Omega + \varepsilon 'e_N)$ that vanishes on $T(3r/4, \rho ) \cap    (\partial \Omega + \varepsilon '\, e_N)$ provided $ \varepsilon '$ is sufficiently small. Since the equation $ \Delta  u-g(x,u)=0$ has its absorption term continuous in $\overline  {\omega _{ \varepsilon '}}\times { \mathbb R} $ a standard result asserts that there is a minimum solution $w_ { \varepsilon ,\varepsilon '}$, with $(v- \varepsilon )_+ \leq w_{  \varepsilon ,\varepsilon '} \leq u$ and that $w_ { \varepsilon ,\varepsilon '}=0$ on
$T(P,r) \cap    (\partial \Omega + \varepsilon '\, e_N)$.

Now by part, a) (and Harnack) we have on $T(r/2, \rho /2) \cap (\Omega + \varepsilon 'e_N)$
$$(v(x)- \varepsilon )_+ \leq w_{ \varepsilon , \varepsilon '}(x) \leq cu(A') r^{N-2} \mathbb G _{A'}^{\varepsilon '}(x) \leq  cu(A') r^{N-2} \mathbb G _{A'}(x) $$
Notice that $c$ here is independent of $ \varepsilon $ and $ \varepsilon '$. So $(v(x)- \varepsilon )_+ \leq cu(A') r^{N-2} \mathbb G _{A'}(x) $ and letting $ \varepsilon $ goes to zero we get that the set of all the considered functions $v$ uniformly goes to zero on $\partial \Omega  \cap T(r/2, \rho /2)$.

c) Finally the supremum of all the nonnegative continuous subsolutions $v$ vanishing on $\partial \Omega \setminus F$ and dominated by $u$ vanishes on $\partial \Omega \setminus F$. It is also well known to be a (continuous) solution and the proposition is proved. $\square$

The next statement complements Proposition \ref{aux2}
\begin{proposition}\label{p:U_F} In addition to the basic conditions, suppose that $g$ satisfies
 \begin{equation}\label{g_growth}   g(x_1,t)/t\to\infty \quad\text{as }t\to\infty, \,  \end{equation}
for some $x_1 \in   \Omega $. Then there is a largest solution $U_F$ of  \eqref{geq} in $\Omega $ such that $U_F\in C(\bar\Omega\setminus F)$ and $U_F=0$ on $F':=\partial\Omega\setminus F$
\end{proposition}
{\sl Proof.} Let ${ \mathcal U}$ denote the set of all solutions of  \eqref{geq} in $\Omega $. Then $\{h\circ v:\,v\in \mathcal{U}\}$ is locally  uniformly bounded in $\Omega $ by  \eqref{Vcond}, so that  by  \eqref{g_growth}, $\sup  \{ v(x_1)\,;\, v \in   \mathcal{U} \}< \infty $. Since  the potentials $V_v:=h \circ v$ are locally uniformly bounded, the solutions $v$ satisfy Harnack inequalities uniformly w.r.\ to $v \in   { \mathcal U}$ in $\Omega $. Therefore $U:=\sup \,{ \mathcal U}$ is locally bounded and in fact a solution of  \eqref{geq}. This proves the proposition for $F=\partial \Omega $. The general case follows using Proposition \ref{aux2} with $u=U_{\partial \Omega} $.

 \section{$(g, \Delta  )$-Moderate solutions and $g$-good measures}
In this section we assume that $g$ satisfies conditions  \eqref{g-basic}  and  \eqref{Vcond}.
Additional assumptions will be mentioned if needed.

\begin{theorem} \label{apgdelta}
Let $u$ be a positive $g$-moderate solution of  \eqref{geq} in the bounded Lipschitz domain $\Omega $. If
$\nu:={\mathrm tr}_{\partial \Omega }\,u$ then there exists a positive measure $\tau \in   \mathfrak M_+(\partial \Omega )$ such that $  \tau \le\nu$, $   \tau \ne 0$  and $\tau\in {\mathfrak M}_+^{g, \Delta  }$, i.e., (recall  $\psi (x)=1\wedge G(x_0,x)$)
\begin{equation}\label{strong1}
 \int_\Omega g\circ \mathbb K [\tau]\;\; \psi\;dx<\infty.
\end{equation}
\end{theorem}

\proof Let $V(x):=h(x,u(x))$, $x \in   \Omega$ .  By assumption \begin{equation}\label{Xs-mod1}
 \int_\Omega   (g\circ u)\;\psi\,dx=\int_\Omega   uV\psi\,dx<\infty.
\end{equation}
and $u$ is  a moderate $L^V$ harmonic function in $\Omega $. Denote by $\nu'\in \mathfrak M_+(\partial \Omega )$  the measure such that $u=\mathbb K^V_{\nu'}$ and recall that  by Proposition \ref{MBLB} the measures $\nu'$ and $\nu$ are mutually absolutely continuous. By  Proposition \ref{A2} and Fubini's theorem,  $\nu '$-almost every point $ \zeta   \in   \partial \Omega $ is $L^{V }$ regular. Thus using Proposition \ref{sing-reg} and  fixing some $ \theta   \in   (0,\frac  \pi  2)$ sufficiently small (depending on $\Omega $)  we may find a compact subset $  F$ of $\partial \Omega $ with $\nu(  F)>0$ and positive constants $c_1$ and $r$ such that whenever $ \zeta   \in     F$ there is a point $x=x( \zeta  ) \in   \Omega $ such that $ \vert   \zeta  -x \vert  =r$,  $( \zeta  ,x]$ is $ \theta  $ non-tangential at $ \zeta  $ for $\Omega $ and
$$\mathbb G (x_0,z) \leq c_1 \mathbb G ^V(z,x_0),\;\;  \quad   \forall   z \in   ( \zeta  ,x  \rangle  .$$

Set $\nu '_{F}=\mathbf 1_{F}\,\nu' $. Then $v:=\mathbb K^V_{\nu '_F}$ is a positive subsolution of (1.1) and the maximal solution $w$ dominated by $u$ and vanishing on $\partial \Omega \setminus   F$ is moderate, its $m$-boundary trace being $ \nu_{F}:=\mathbf 1_{F} \nu $. This follows from Lemma \ref{KVK} which implies that $\mathbb K^V_{\nu '_{F}} \leq w \leq \mathbb K_{\nu_{F}}$ in $\Omega $ and then that ${\mathrm {tr}}_{\partial \Omega }(w)=\nu_{F}$ .

In particular if $W(x):= h(x,w(x))$,  $x \in   \Omega $, we have $w=\mathbb K^{W}_{\lambda}$ for a positive measure $ \lambda  \in   \mathfrak M _+(\partial \Omega )$ which is equivalent to $\nu_F$ and hence to $\mathbf 1_{F} \nu'$. (The measures $\lambda,\nu_F$ are equivalent in the sense that they are mutually absolutely continuous.)

Also, since $0 \leq W \leq V$, we again have $$\mathbb G (x_0,z) \leq c_1 \mathbb G ^W(x_0,z)$$
whenever $z \in   ( \zeta  ,x]$ with $ \zeta   \in   F$ and $x=x( \zeta  )$ as above. By Proposition \ref{FMG} this means that for some other constant $c_2 \geq 1$
$$(c_2)^{{-1} }\mathbb K_{ \zeta  }^W(z)  \leq \mathbb K_{ \zeta  }(z) \leq c_2 \,\mathbb K_{ \zeta  }^W(z),  \quad    \quad   \forall   z \in   ( \zeta  ,x (\zeta  )], \;  \zeta   \in    F$$
and by Proposition \ref{mainSec4}  above --notice that $W \leq c\, [d(x, F)]^{-2}$-- we conclude that for some other constant $c_3>1$ we even have
$$(c_3)^{{-1} }\mathbb K_{ \zeta  }^W(z)  \leq \mathbb K_{ \zeta  }(z) \leq c_3 \; \mathbb K_{ \zeta  }^W(z),  \quad    \quad   {\mathrm {\ for \ all \ }}  \zeta   \in  F      \mathrm{\ and \ all  } \; z \in   \Omega  $$
So that finally, taking $  \tau \not \equiv 0$ with  $ c_3\tau  \leq  \nu' \wedge \lambda $
$$\int (g \circ \mathbb K _ {\tau} )\,  \psi \, dx \leq \int (g \circ \mathbb K ^{W} _{c_3\tau})  \,  \psi \, dx
 \leq \int (g \circ \mathbb K ^{W} _{\lambda})  \,  \psi \, dx  \leq  \int (g \circ u)  \,  \psi \, dx < \infty .  \quad  \square$$

Based on Theorem \ref{apgdelta} above, the next result  provides a complete characterization of  conditionally g-removable sets. In addition, if $g$ satisfies the $\Delta_2$ condition, we obtain a full characterization of   $g$-good measures and show that every $g$ moderate solution is the limit of an  increasing sequence of $(g,\Delta)$ moderate solutions. In the absence of the $\Delta_2$ condition we establish a modified version of the last two results

At this point it is convenient to introduce the following notation:
\begin{equation}\label{CT}
    {\mathcal T}^{g,\Delta}:=\cup_1^\infty n{\mathfrak M}^{g,\Delta}_+(\partial\Omega).
\end{equation}
Thus \ ${\mathcal T}^{g,\Delta}$ is the convex cone generated by ${\mathfrak M}^{g,\Delta}_+(\partial\Omega)$. If $\mu_1,\ldots,\mu_n\in {\mathfrak M}^{g,\Delta}_+$ then
$$\frac{1}{n}\sum_1^n\mu_i\in {\mathfrak M}^{g,\Delta}_+$$
 and ${\mathcal T}^{g,\Delta}$ may be described as the family of all finite sums of measures from ${\mathfrak M}^{g,\Delta}_+$. Note that if $\nu_1,\nu_2\in {\mathcal T}^{g,\Delta}$ and $\nu_1\leq \nu_2$ then $\nu_2-\nu_1\in {\mathcal T}^{g,\Delta}$.

Clearly, if $g$ satisfies the $\Delta_2$ condition then
${\mathcal T}^{g,\Delta}={\mathfrak M}^{g,\Delta}_+.$

\begin{theorem}\label{gmod1}   A compact set $F  \subset    \partial \Omega $ is  conditionally $g$-removable  if and only if  $F$ is ${\mathfrak M}^{g,\Delta}_+$ null.
\end{theorem}
\proof If $F$ is not ${\mathfrak M}^{g,\Delta}_+$ null then, by definition, there exists a positive measure $\tau\in    {\mathfrak M}_+^{g,   \Delta  }$ which is concentrated on $F$. Thus $ \tau $ is $g$-good and $   \mathbb S^g_\tau$ is a positive moderate solution of  \eqref{geq} that vanishes on $\partial \Omega \setminus  F$.

Suppose now that $F$ is ${\mathfrak M}^{g,\Delta}_+$ null. Let $u$ be a non-negative $g$-moderate solution that vanishes on $\partial \Omega \setminus  F$
and let $\mu:=  {\mathrm tr}\,u$. Clearly $\mu(\partial \Omega \setminus F)=0$.
If $u$ is positive then, by Theorem \ref{apgdelta}, there exists a positive measure $\tau\in {\mathfrak M}^{g,\Delta}_+$ such that $\tau\leq \mu$. But this is impossible because $F$ is ${\mathfrak M}^{g,\Delta}_+$ null.  Thus $u=0$.
 \qed

\begin {theorem}\label{gmod2}
If $u$ is a positive $g$-moderate solution of  \eqref{geq} in $\Omega $ then $u$ is the limit of an increasing sequence  of solutions of  \eqref{geq} whose boundary trace belongs to ${\mathcal T}^{g,\Delta}$. Consequently, there exists a sequence  $\{\mu_n\}\subset{\mathfrak M}^{g,\Delta}_+$  such that
\begin{equation}\label{CTmod}
    {\mathrm tr}(u)=\sum_1^\infty\mu_n
\end{equation}

If $g$ satisfies the $\Delta_2$ condition then every
 $g$-moderate solution is the limit of an increasing sequence of
$(g, \Delta  )$~moderate solutions.
\end{theorem}

\proof Let $u$ be a positive $g$-moderate solution with m-boundary trace  $\mu$.  Let $u^*$ be the supremum of all $g$-moderate solutions $v$ such that ${\mathrm tr} (v)\in {\mathcal T}^{g,\Delta}$ and $v\leq u$. Clearly $u^\ast $ is $g$-moderate (see Proposition \ref{gmodstab}), $\mu  ^\ast:={\mathrm tr}(u^*) \leq \mu  $ and $\mu  ^*$ is the upper envelope of all $\nu  \in   {\mathcal T}^{g,   \Delta  }$ which are majorized by  $\mu  $.

Suppose that  $u^*\not\equiv u$ so that $\mu^*:=  {\mathrm tr}\, u^*\not \equiv \mu$. Then $\mu-\mu^*$ is a nonzero  $g$-good measure and by     Theorem \ref{apgdelta}  there exists a nonzero  measure $\tau'\le\mu-\mu^*$    in  $    {\mathfrak M}_+^{g,   \Delta  }$. If  $\nu\in {\mathcal T}^{g,   \Delta  }$ then $\nu+\tau'\in {\mathcal T}^{g,   \Delta  }$. In addition $\mu^*+\tau'\leq \mu$. Therefore
$\mu^*+\tau'\leq \mu^*$, which contradicts the fact that $\tau'$ is positive.
 Hence, $u^*=u$. Finally, every positive measure $\mu$ that is the limit of an increasing sequence of measures in ${\mathcal T}^{g,\Delta}$ can be represented as in  \eqref{CTmod}.

Now assume that $g$ satisfies the $\Delta_2$ condition. In this case ${\mathfrak M}^{g, \Delta  }_+={\mathcal T}^{g, \Delta  }$. Therefore
the last assertion of the theorem is a consequence of the first.
 \qed

\begin{theorem}\label{gmod3}
(i) If $\mu\in    {\mathfrak M}_+$ is a $g$-good measure then $\mu(A)=0$ for every ${\mathfrak M}^{g,\Delta}_+$ null set $A$.

\noindent(ii) If $\mu\in    {\mathfrak M}_+$ vanishes on ${\mathfrak M}^{g,\Delta}_+$ null sets then there exists a  Borel function $f:\partial\Omega\mapsto(0,1]$  such that $f\mu$ is a $g$-good measure.

\noindent(iii) If $g$ satisfies the $\Delta_2$ condition then, $\mu\in    {\mathfrak M}_+$ is a $g$-good measure if and only if  $\mu(A)=0$ for every ${\mathfrak M}^{g,\Delta}_+$ null set $A$.
 \end{theorem}

\proof (i) Let $A$ be a compact ${\mathfrak M}^{g,\Delta}_+$ null set. If $\mu$ is $g$-good and $\mu(A)>0$ then $\tau:=\mu\mathbf{1}_A$ is a $g$-good positive measure. Since $A$ is compact, ${\mathbb S}^g[\tau]$ is a positive g-moderate solution that vanishes on $\partial\Omega\setminus A$. By Theorem \ref{gmod1} $u=0$, which is absurd. This implies (i).

(ii) If $\mu\in    {\mathfrak M}_+$ vanishes on ${\mathfrak M}^{g,\Delta}_+$ null sets then, by    Theorem \ref{rich}, it is the limit of an increasing  sequence of measures belonging to $  { {\mathcal T}}_+ ^{g, \Delta  }(\partial  \Omega )$ (see Corollary \ref{richExample}). Therefore it can be represented by a series $\sum_1^\infty\mu_n$ where $\mu_n\in {\mathfrak M}^{g,\Delta}_+$. Every measure in ${\mathfrak M}^{g,\Delta}_+$ is $g$-good and the maximum of two $g$-good measures is $g$ good (see Proposition \ref{gmodstab}).

Therefore
$\nu_n:=\max(\mu_1,\cdots,\mu_n)$ is $g$-good. Consequently
$\nu:=\lim_{n\to\infty}\max(\mu_1,\cdots,\mu_n)$
is $g$ good, $\nu\leq \mu$ and $\nu$ is absolutely continuous w.r. to $\mu$.

(iii) Now assume that $g$ satisfies the $\Delta_2$ condition. Recall that in this case ${\mathfrak M}^{g, \Delta  }_+={\mathcal T}^{g, \Delta  }$.
If $\mu\in    {\mathfrak M}_+(\partial \Omega  )$ is a  positive measure that vanishes on $  { \mathfrak M}_+ ^{g, \Delta  }$-null sets then, by    Corollary \ref{richExample2}, it is the limit of an increasing  sequence of measures belonging to $  { \mathfrak M}_+ ^{g, \Delta  }(\partial  \Omega )$. Since every measure in $  { \mathfrak M} _+^{g,   \Delta  }(\partial \Omega  )$ is $g$-good it follows that $\mu$ is $g$-good.\qed
 \quad

\section{Positive solutions and $g$-moderate solutions}
As in the previous section it will always be assumed that $g$ satisfies conditions  \eqref{g-basic} and  \eqref{Vcond}.

Let $u$ be a positive solution of  \eqref{geq} in $\Omega $. Let $V$ be defined as in  \eqref{Vdef} and let $\mu$ be the $L^V$ boundary measure for $u$, i.e., $u= { \mathbb  K} ^V_\mu$.

 Set
\begin{equation} \label{u*1}u^\ast=\sup  \{ v\,;\, v \mathrm{ \ is \  a \ g-moderate\ solution\ in \ }\Omega \mathrm {\  and \ } \,0 \leq v \leq u\,  \}.\end{equation}
The function $u^\ast $ is  the largest $ \sigma  $-moderate solution of  \eqref{geq} dominated by $u$.

Before stating our final result we recall a classical lemma on decomposition of measures that is used in our proof.

\begin{lemma}\label{mu-cap} (\cite{FTS}, see also \cite[Appendix 4A]{BMP04})  Let $\mu$ be a finite measure on a measurable space $(Y,{ \mathcal F})$
and let ${ \mathcal G }$ be a subset of ${ \mathcal F}$  such that :\\[1mm]
1. $ \emptyset   \in   { \mathcal G}$ and  ${ \mathcal G }$ is closed  with respect to  finite or countable unions,\\[1mm]
2. $A\in { \mathcal G }$, $A' \in   { \mathcal F}$ and $A'\subset A   \Longrightarrow  A'\in { \mathcal G}$.\\[1mm]
\quad \hspace{4mm}    Then $\mu$ can be uniquely written in the
 form
$  \mu=\mu_0+\mu_1, $ where $\mu_0$ and $\mu_1$ are finite measures on $(Y,{ \mathcal F})$,  $\mu_1(A)=0 $    for all $    A\in { \mathcal G}$ and $ \mu_0$  is \ concentrated  on  a \ set $A_0\in { \mathcal G}.$

\end{lemma}

Combining this lemma and Theorem \ref{rich} we obtain,

   \begin{proposition}\label{split-mu} Let $\tau\in \mathfrak M_+( \partial \Omega )$ be a positive finite  Borel measure. Then:

(i) There is a unique decomposition $\tau=\tau_d+\tau_c$
 with $\tau_d,\,  \tau _c\in   \mathfrak M_+(\partial \Omega )$, $ \tau _d$  diffuse and $\tau_c$  singular with respect to  $C_{g,  \Delta }$.

(ii) $ \tau _d=\sum_1^\infty  \tau _n$ for some increasing   sequence $\{\tau_n\}$   in $\mathfrak M_+^{g,\Delta}(\partial \Omega )$
\end{proposition}
\proof The first statement  is a special case  of  Lemma \ref{mu-cap}.  The second statement follows from  Theorem \ref{rich} (see Corollary \ref{richExample}) by the same argument as in the proof of Theorem \ref{gmod2}.
   \qed

If $F$ is a closed subset of $\partial \Omega $, we denote
 (recall $u$ is a nonnegative solution of  \eqref{geq} in $\Omega $),
\begin{equation}\label{[u]F}
   [u]_F:=\sup\{v\in \mathcal{U}_F:\, v\leq u\}.
\end{equation}
where ${ \mathcal U}_F$ is the set of all non-negative solutions of  \eqref{geq} in $\Omega $ vanishing on $\partial \Omega \setminus F$. By Proposition \ref{aux2}, $[u]_F$ is the largest solution dominated by $u$ and vanishing on $\partial \Omega \setminus F$.

We now state this section's main result. Recall that $u$ is a positive solution of ($\ref{geq}$) and that $u^\ast$ is defined by (\ref{u*1}).

\begin{theorem}  \label{gsmod}  If $\mu$ is not singular  relative to $C_{g,\Delta}$ then $u$ dominates a positive $g$ moderate solution. Furthermore, there exists a Borel function $f:\partial \Omega  \to  (0,1]$ such that $f\mu_d$ is $g$-good and
  \begin{equation} \label{u*3}
    {\mathbb  S} ^g_{f\mu  _d}\leq u.
    \end{equation}

Let $u^\ast$ be as in  \eqref{u*1}. For every compact subset $F$ of ${ \mathcal R}(u):={ \mathcal R }^V$,
\begin{equation} \label{u*4}
[u]_F=[u^\ast]_F.
\end{equation}
\end{theorem}

With the next lemma we begin the proof of Theorem \ref{gsmod}. This lemma is an adaptation of \cite[Lemma 5.5]{MM-mod}.  For the convenience of the reader we reproduce the proof here.

  \begin{lemma}  Let $\nu\in {\mathfrak M}^{g,  \Delta }( \partial \Omega )$ be a positive measure.  Suppose that there exists no positive solution of  \eqref{geq} dominated by the supersolution $v=\inf(u, {  \mathbb  K} _\nu)$. Then $\mu\bot\nu$.
\end{lemma}

\proof  Set $V'=h(v)$. Clearly $v$ is  $L^{V'}$ superharmonic since  $v$ is a supersolution of  \eqref{geq}.

 Moreover \emph{ $v$ is an $L^{V'}$ potential.} For if $w$ is a  nonnegative $L^{V'}$ harmonic function dominated by $v$, then $w$ is a subsolution of  \eqref{geq}. Hence, since $v$ is a supersolution,  there exists  a solution $w'$ of  \eqref{geq} with $w \leq w' \leq v$. By the assumption $w'=0$ and hence $w=0$. This shows that $v$ is an $L^{V'}$ potential.\medskip

 Evidently, the harmonic function
$ { \mathbb  K} _\nu$ is  $L^{V'}$  superharmonic. Since
$$\int_{\Omega   } {  \mathbb  K} _\nu V'    \psi  \, dx\leq \int_{\Omega   }g( {  \mathbb  K} _\nu )   \psi  \,dx<\infty.$$
$ {  \mathbb  K} _\nu$ is moreover an   $L^{V'}$ - moderate  nonnegative superharmonic function.  By   Lemma \ref {A0}, the largest $L^{V'}$ harmonic dominated by $ {  \mathbb  K} _\nu$, say $w$,  is $L^{V'}$ moderate and has the m-boundary trace $\nu$. Of course
$$p:=  {  \mathbb  K} _\nu-w$$
is an $L^{V'}$-potential and by the results in \cite{Annals87} (see Section 3 above)
 $$w= {  \mathbb  K} ^{V'}_{\nu'}$$ where $\nu'$ is a positive finite measure on $ \partial \Omega $. By     Proposition \ref {MBLB}, $\nu$, $\nu'$ are mutually absolutely continuous.

 By  the relative Fatou theorem (for $L^{V'}$), since $v,p$ are $L^{V'}$ potentials and $w$ is  $L^{V'}$ harmonic,
\begin{equation}\label{fine-nu}
 v/w\to 0,    \quad    {  \mathbb  K} _\nu/w\to 1     \quad  \text{$L^{V'}$-finely $\nu'$-a.e.}\,.
\end{equation}
Because $v=\inf(u, {  \mathbb  K} _\nu)$,  \eqref{fine-nu} implies that
\begin{equation}\label{fine-nu'}
  u/w \to 0    \quad  \text{$L^{V'}$- finely $\nu'$-a.e.}
\end{equation}
Further, by   \eqref{fine-nu} and  \eqref{fine-nu'}
\begin{equation}\label{fine-u}
 \frac{u}{ {  \mathbb  K} _\nu} \to 0    \quad  \text{$L^{V'}$- finely $\nu'$-a.e.}
\end{equation}
Thus, since  the function $ u/ {  \mathbb  K} _\nu$ satisfies the strong Harnack condition in $\Omega $,
\begin{equation}\label{nt-nu}
 \frac{u}{ {  \mathbb  K} _\nu} \to 0   \quad  \text{n.t. $\nu$-a.e.}
\end{equation}
where we have also used the fact that   $\nu$ and $\nu'$ are equivalent.

So ${  \mathbb  K} _\nu/u \to  + \infty $ n.t. at $\nu$-a.e.\ point in $\partial \Omega $. But on the other hand, viewing $ {  \mathbb  K} _\nu$ as a non negative $L^V$-superharmonic function, the $L^V$ relative Fatou-Doob theorem says that $ {  \mathbb  K} _\nu/u$ has a finite n.t. limit $\mu  $-a.e.\ in $\partial \Omega $. These two asymptotic behaviors for ${  \mathbb  K} _\nu/u$ imply that $\mu \bot \nu$. $\square$

Recall that by definition $u=\mathbb K_\mu ^V$ but that $\mu  $ is not necessarily $g$-good.

\begin{corollary}\label{nu-ac-mu} Suppose that the measure $\nu \in   {\mathfrak M}_+^{g,  \Delta }( \partial \Omega )$ is absolutely continuous with respect to $\mu  $. Then there exists a Borel function $f$ on $ \partial \Omega $  such that $0<f$ $\nu$-a.e. and
\begin{equation}\label{fnu}
    v:= {  \mathbb  S} ^g_{f\nu}\leq u.\end{equation}
    Moreover when $u $ admits a trace $\mu'$ we can choose $f$ such that $f\nu=\nu \wedge \mu  '$.
\end{corollary}
{\sl Proof.}  By the previous Lemma, for every positive measure $\nu' \leq \nu $, there is a positive solution $w$ of  \eqref{geq} in $\Omega $ such that  $0<w \leq \mathbb K_{\nu'}\wedge u$. It follows that the  least upper bound of the set of good measures $ \lambda  \leq \nu$  such that $\mathbb S_ \lambda^g  \leq u$ is a $g$-good measure in the form $f\nu$ with $0<f \leq 1$ $\nu$-a.e.
\qed

\medskip
 \begin{lemma}\label{regmod} If $\mu  $-almost every  $ \zeta   \in   \partial \Omega $ is $L^V$- regular, then $u$ is $ \sigma  $-moderate.
 \end{lemma}

\proof  We exhaust $\mu  $  with an increasing  sequence $ \{ F_n \}$ of compact subsets of $\partial \Omega $ chosen as follows: $F_n$ is the set of all $ \zeta   \in   \partial \Omega $ such that $ {  \mathbb  K} _ \zeta  ^V \leq n  {  \mathbb  K} _ \zeta  $ in $\Omega $.
 The fact that $\mu  (\partial \Omega \setminus F_n) \to  0$ follows from the assumptions and Definition \ref{V-regular}.

Let $\mu_n=\mu\mathbf{1}_{F_n}$ and put $u_n= {  \mathbb  K} _{\mu  _n}^V$.
Then $u_n\uparrow u$ and $u_n$ is a subsolution of  \eqref{geq} dominated by $n\mathbb K_\mu$. Let $\hat u_n$ be the smallest solution of  \eqref{geq} larger than $u_n$.
Then $\hat u_n< n\mathbb K_\mu$ and consequently it is $g$-moderate.
Obviously $\hat u_n\uparrow u$. \qed

    \begin{corollary} \label{FinR} If $F \subset   \partial \Omega $ is compact and  $F   \subset    { \mathcal R}(u)$ then $[u]_F$ is  $\sigma  $-moderate.
 \end{corollary}
\proof We may assume that $v=[u]_F>0$. Consider $V_F=h\circ v$. As $V_F \leq V$ it follows that $\mathbb G ^V \leq \mathbb G ^{V_F}$. Therefore, by Proposition \ref{sing-reg} (iii), every  $L^V$ regular boundary point is also  $L^{V_F}$ regular. Thus $L^{V_F}v=0$, $v$ and $V_F$ vanish on $ \partial \Omega    \setminus  F$ and every $   \zeta  \in  \partial \Omega $ is $L^{V_F}$ regular.
Therefore the result is a consequence of the previous Lemma. \qed

\noindent{\sl Proof of   Theorem \ref{gsmod}.}
By     Theorem \ref{gmod3} (ii) there exists a  Borel function $f_0:\partial\Omega\mapsto(0,1]$  such that $f_0\mu_d$ is a $g$-good measure. By     Theorem \ref{gmod2} there exists a sequence $\{\tau_n\}\subset {\mathfrak M}^{g,\Delta}_+$ such that
\begin{equation}\label{fmud}
f_0\mu_d=\sum_1^\infty\tau_n.
\end{equation}
By  Corollary \ref{nu-ac-mu}, for every $n\in  {  \mathbb  N} $ there exists a function $h_n: \partial \Omega \mapsto (0,1]$  such that $h_n:\partial\Omega\mapsto (0,1]$ and
\begin{equation*}
    v_n:= {  \mathbb  S} ^g_{h_n\tau_n}\leq u.
\end{equation*}
Further, $\nu_n:=\max(h_1\tau_1,\cdots,h_n\tau_n)$  is $g$-good, $\nu_n\leq \mu_d$ and
 $$w_n={  \mathbb  S} ^g_{\nu_n}\leq u.$$
 Therefore $\nu$ is $g$-good and
 $$w:=\lim w_n = {  \mathbb  S} ^g_{\nu}\leq u.$$
Finally  \eqref{fmud} implies that $\mu_d$ is absolutely continuous w.r. to $\nu$. This proves  \eqref{u*3}

By     Corollary \ref {FinR}, if $F$ is a compact subset of $ { \mathcal R}(u)$ then $[u]_F$ is $ \sigma  $-moderate . Therefore
$[u]_F\leq [u^*]_F.$
Since $u^*\leq u$ we obtain  \eqref{u*4}. \qed

\newpage

 \renewcommand\thesection{\Alph{section}}
\setcounter{section}{1}

\noindent \hspace{6mm}
\noindent \hspace{6mm}

      \vspace{-2truemm}

 \section*{{\bf Appendix A. Some classes of ``good" absorption terms
  } \\ by Alano Ancona }

 \renewcommand\thesection{\Alph{section}}
\setcounter{section}{1}
\setcounter{equation}{0}
\setcounter{theorem}{0}

\noindent \hspace{6mm}
\noindent \hspace{6mm}

      \vspace{2truemm}

{\textsc {Summary.}} {\small  We first describe in section \ref{A1} a class of absorption terms $g(x,t)$ such that the condition  \eqref{VcondSp} from section 1 holds (i.e.\ every positive solution $u_0$ of  \eqref{theequation} below is a solution of a Schr\"odinger equation that behaves well with respect to dilations). Classes of functions $g$ such that the stronger condition  \eqref{Vcond} holds is later described in section \ref{A22}.  In section  \ref{A44}   Theorem \ref{thap2}  gives (for a somewhat more restricted class of $g$) a simpler way to condition    \eqref{Vcond} by deriving it from condition   \eqref{VcondSp}.

Recall that condition  \eqref{VcondSp} implies both (uniform inner) Harnack inequalities for $u$ as well as a useful uniform boundary Harnack principle. See remark \ref{schrrem}.

    The exposition is made essentially independent of the previous sections.}

  \vspace{2mm}

  \subsection{}\label{A1}
    We consider in a domain $\Omega $ of ${ \mathbb  R}^N$, $N \geq 1$, a semi-linear equation
\begin{equation}\label{theequation}-\Delta  u(x)+g(x,u(x))=0\end{equation}

where the  function $g:  \Omega \times { \mathbb R}_+  \to  { \mathbb R}_+ $ is  continuous.

 A nonnegative  function $u$ in $\Omega $ is a solution of  \eqref{theequation} if $u$ and $g(.,u(.))$ are in $L^1_{\rm loc}(\Omega )$  and  $ \Delta  u=g(.,u(.))$ --in the distribution sense--. In particular $u$ is subharmonic -thus locally bounded- and $ \Delta  u \in   L^ \infty _{{\rm loc}}(\Omega )$. So $u$ is continuous and  in $ W^{2,p}_{{\rm loc}}(\Omega )$ for all $p< \infty $.

It is once for all  assumed that $g$ satisfies the following three assumptions.

{\bf ($ND$)} \  for each $x \in   \Omega $, $t \mapsto g(x,t)$ is $ \geq 0$ nondecreasing  on $[0,+ \infty )$,

{\bf  ($H$)}  there is a real $C_1 \geq 1$ such that $g(x',a) \leq C_1\,  g(x,a)$ whenever $ \vert  x'-x \vert   \leq {\frac  { 1} {2}} \delta (x)$.

Here $ \delta (x):=\mathrm{dist}(x,{ \mathbb  R}^N\setminus \Omega )$, and ($H$)  is reminiscent of Harnack inequalities.  It is also assumed that  there are positive constants $C_2$, $c_0$ such that the following holds:

{\bf ($KOT$)} If $x \in   \Omega $, $a \in   { \mathbb R} _+^*$ are such that  ${\frac  { g(x,a)} {a}} > c_0\,  \delta (x)^{-2}$, then
\begin{equation}\label{cond1}\int _{2a}^ \infty  {\frac  { dt} {\sqrt{G(x,t)} } } \leq  C_2\, \sqrt{  {\frac  { a} {g(x,a)}} },\end{equation}
where $G(x,t)=\int _0^t\, g(x,s)\, ds$.  When $g$ is $x$-independent the finiteness -for some $a>0$- of the  l.\ h.\  s.\ of  (\ref{cond1}) is  the well-known Keller-Osserman condition (see the comments  in the beginning of section 1).

Note that ($KOT$)  forces $g(x,0)=\lim_{t>0,t   \to      0}g(x,t)=0$.

 \begin{proposition}\label{ConcCond-St} Under the above conditions {\bf ($ND$)}, {\bf  ($H$)} and {\bf ($KOT$)}, there is a positive constant  $\bar c$ such that for every nonnegative solution $u$ in $\Omega $ of equation (\ref{theequation})
\begin{equation}\label{Potcond}g(x,u(x)) \leq \bar c \, u(x) \, \delta (x)^{-2} \quad \mathrm{ for \ all\ } x \in   \Omega \end{equation}
The constant $\bar c$ may be chosen depending only on $N$, $c_0$, $C_1$ and $C_2$.
 \end{proposition}

 We first note the following elementary fact.

 \begin{lemma} \label{lemelem} Set $ C_3=2+C_2$. Under the above assumptions,  we have  for $x \in   \Omega $ and $a \in   { \mathbb R} _+^*$ such that  ${\frac  { g(x,a)} {a}}> c_0\,  \delta (x)^{-2}$,
 \begin{equation}\label{cond2} \int _{a}^ \infty  {\frac  { dt} {\sqrt{G(x,t)-G(x,a)} } } \leq  C_3\, \sqrt{  {\frac  { a} {g(x,a)}} }.\end{equation}
 \end{lemma}

 {\sl Proof.} By the convexity of $t \mapsto G(x,t)$ in ${ \mathbb R} _+$, $G(x,t)-G(x,a) \geq g(x,a)(t-a)$. Thus $$\int _{a}^ {2a}  {\frac  { dt} {\sqrt{G(x,t)-G(x,a)} } } \leq 2\,\sqrt{  {\frac  { a} {g(x,a)}} }.$$
On the other hand for $t \geq 2a$,  $G(x,t) \geq 2\,G(x,a) $ since $G(x,t)/t$ is nondecreasing in $t$. Whence $G(x,t)-G(x,a) \geq {\frac {1}Ê{2} }G(x,t)\geq  G(x,a)$ and
$$\int _{2a}^ \infty  {\frac  { dt} {\sqrt{G(x,t)-G(x,a)} } } \leq \sqrt 2 \int _{2a}^ \infty  {\frac  { dt} {\sqrt{G(x,t)} } } $$
On using (\ref{cond1}) the result follows. \qed

 \bigskip

  The next lemma appears in   \cite{DDGR} -under some slightly stronger assumptions -. Consider an ordinary differential equation
   \begin{align}  y''(r)+{\frac  { N-1} {r}}\,y'(r)&=f(y(r))
 \tag{*}\nonumber
\end{align}
 with  $f \geq 0$ and non decreasing in  ${ \mathbb R} _+$. Set $F(t)=\int _0^t f(s)\, ds$.

 Recall that a nonnegative subsolution of $(*)$ in $(a,b)$, $0 \leq a<b<+ \infty $, is a nonnegative locally bounded function $y$ in $(a,b)$  such that  the distribution $[r^{N-1}y'(r) ]'$ is a positive measure in $(a,b)$ larger than $r^{N-1}f(y(r))\, dr$. So $y$ may be viewed as locally Lipschitz function in $(a,b)$.

\begin{lemma} \label{4Lemma} {\rm (See \cite[Lemma 2.7]{DDGR} )}  Suppose that $y$ is a  nondecreasing and nonnegative subsolution of  $(*)$  in a neighborhood of $[r_0,r_1]$, $0 < r_0 < r_1< \infty $. Suppose also that $F(y(r_0))>0$. Then,
\begin{equation}\label{4lemma} \int _{y(r_0)}^{y(r_1)}{\frac  { du} {\sqrt {2(F(u)-F(y(r_0))} } }  \geq {\frac  { \,r_0} {N-2}}\, [1- {\frac  { r_0^{N-2}} {r_1^{N-2}}}] \; ({\rm  or \ }  \geq r_0 \log {\frac  { r_1} {r_0}}, {\rm \ if \ } N=2).\end{equation}
\end{lemma}

To prove the lemma one may easily  adapt the argument given in \cite{DDGR} for smooth $y$. We skip the details.

\medskip

\begin{corollary} \label{Cor4}If there exists  a nonnegative  and nondecreasing subsolution $y$ in  $(0,T)$ with $F(y(0+))>0$, then
\begin{equation}\label{4lemmaCons} \;{\frac  { 2^{N-2}-1} {(N-2)2^{N-2}}}\; T \leq \sqrt 2\int _{y(0)}^{ \infty }  {\frac  { du} {\sqrt {F(u)-F(y(0))} } } \end{equation}
where the l.h.s.\ should be replaced by $ ( \log 2)\,T\,$ when $N=2$.
\end{corollary}

\medskip  The corollary  follows from the the lemma by taking $r_0=T/2$, $r_0<r_1<T$, and letting $r _1\to  T$.
(Notice that $ \int _a^ {+ \infty } {\frac  { dt} {\sqrt {F(t)-F(a)}}}$ is a decreasing function of $a$, since $F$ is convex.)

\medskip

 {\sl Proof of Proposition \ref{ConcCond-St}.} Consider a point $x_0 \in   \Omega $ and  the ball $B=B(x_0, \delta(x_0) /2)$. Then, by ($H$), $u$ is a nonnegative subsolution of the equation
\begin{equation} \label{eqg(x0)}\Delta  u(x)-{\frac  { 1} {C_1}}  \, g(x_0,u(x))=0\end{equation}
 in $B$. Letting  $g_0(t)={\frac  { 1} {C_1}} g(x_0,t)$ we may rewrite this equation: $ \Delta  u-g_0(u)=0$.

 For every isometry $R$ of ${ \mathbb R} ^N$ that fixes $x_0$, $u \circ R$ is again a subsolution of (\ref{eqg(x0)}) in $B$ and  the function  $w:=\sup  \{ u \circ R\,;\, R$ isometry  fixing $x_0\, \}$ is also a (continuous) subsolution of (\ref{eqg(x0)}) in $B$  in the form $w(x)= \tilde u(d(x,x_0))$ for some continuous $ \tilde u:[0, \delta (x_0)/2] \to  { \mathbb R} _+$. The function $w$ is a classical subharmonic function in $B$ and is invariant under every rotation that fixes $x_0$. Thus $ \tilde u$ is nondecreasing and evidently  $ \tilde u$ is a subsolution of the ODE  : $y''(r)+{ \frac {  N-1} {r }}y'(r)-g_0(y(r))=0$ in $(0, \delta (x_0)/2)$.

 Applying Corollary \ref{Cor4}, we see that if $g_0(u(x_0))>
0$ -hence $u(x_0)>0$- and if say $N \geq 3$, \medskip
 \begin{align} \nonumber {\frac  {  \delta (x_0) } {N-2}}   {\frac  { 2^{N-2}-1} {2^{N-1}}  }   &  \leq \sqrt 2  \int _{u(x_0)}^{ \infty }  {\frac  { du} {\sqrt {G_0(u)-G_0(u(x_0))} } } \\
 \\ &=\sqrt {2 C_1}  \int _{u(x_0)}^{ \infty }  {\frac  { du} {\sqrt {G(x_0,u)-G(x_0,u(x_0))} } }.\nonumber
 \end{align}
 So if  moreover ${\frac  { g(x_0,u(x_0))} {u(x_0)}}  > c_0\,  \delta (x_0)^{-2}$  it follows from  (\ref{cond2}) that
 $$c_1  \delta (x_0) \leq C_3\sqrt{2 C_1} \sqrt{  {\frac  { u(x_0)} {g(x_0,u(x_0))}} }$$
with $c_1= {\frac  { 1} {N-2}}   {\frac  { 2^{N-2}-1} {2^{N-1}}  }$. This implies that  $g(x_0,u(x_0))   \leq c  \delta (x_0)  ^{-2}(x_0)\,u(x_0)$, $c={ \frac {  2C_1(C_3)^2} { c_1\,^2}}$.

The only other possibility is  that  $g(x_0,u(x_0)) \leq c_0 u(x_0) \,  \delta (x_0)^{-2}$.
Proposition  \ref{ConcCond-St} is proved.   \qed

 \begin{remark} \label{schrrem} {\rm Proposition \ref{ConcCond-St} says that $u$ is a solution of a Schr\" odinger equation $L^V(u):= \Delta  u-Vu=0$ for a potential $V \in   L^ \infty _{loc}(\Omega )$, $V(x) \leq \bar c  \delta (x)^{-2}$, $\bar c$ depending only on the constants $C_1$, $C_2$, $c_0$ and $N$. In particular there is a constant $c$ depending only on $\bar c$ such that $u(x) \leq c\, u(x')$ whenever $ \vert  x'-x \vert   \leq {\frac  { 1} {2}}  \delta (x)$ (Harnack inequalities). Moreover, the Harnack boundary inequalities of \cite{Annals87} (see \cite{An-MZ}) are available for $L^V$.}
 \end{remark}

 \begin{remark} \label{ddgrrem} {\rm We note here that for $g$  space independent, the paper  \cite{DDGR} about blow up boundary solutions removes the classical monotonicity assumption on $g$. It is in particular shown there that $g=g(u)$ satisfies the Keller-Osserman condition  (resp. a sharpened Keller-Osserman condition) if and only if there exists a positive boundary blow-up solution on some ball (resp.\ on any smooth bounded domain). See \cite{DDGR}.}
 \end{remark}

\subsection {}  \label{A22}
 We next show that if $\Omega $ is Lipschitz and if $g$ fulfills a supplementary assumption then the conclusion in Proposition \ref{ConcCond-St} can be strengthened : property  \eqref{Vcond} holds, i.e. there is a positive real $\bar c>0$  such that for every closed subset $F$ of $\partial \Omega $ and every nonnegative solution $u$ of (\ref{theequation}) in $\Omega $ that vanishes on $\partial \Omega \setminus F$,
\begin{equation} \label{Potcond2bis} g(x,u(x)) \leq { \frac {  \bar c} { d(x,F)^2}}\, u(x),\;\; x \in   \Omega. \end{equation}

From now on, $\Omega $ is assumed to be bounded and Lipschitz. For $ \varepsilon >0$, define $\Omega ( \varepsilon ):= \{ x \in   \Omega \,;\,  \delta (x)< \varepsilon \,  \}$.

  As before it is assumed that $g$ is continuous in $\Omega \times { \mathbb R} _+$, $g(x,0)=0$  and  that $g(x,.)  $ is nondecreasing for all $x \in   \Omega $; moreover the Harnack type condition ($H$) holds. In other words,  there is a real $C_1 \geq 1$ such that $g(x',a) \leq C_1\,  g(x,a)$ whenever $ \vert  x'-x \vert   \leq {\frac  { 1} {2}} \delta (x)$,

  \begin{definition} We say that  $g$ is quasi inwardly increasing (near the boundary) in $\Omega $ if there are constants $ \varepsilon _0>0$ and $C \geq 1$ such that :

\vspace{-2truemm}  ($QI$) whenever $x \in   \Omega( \varepsilon _0) $,  $0<r< \varepsilon _0$ and $a>0$, we have $g(x,a) \leq C\, g(x',Ca)$ for at least one  point $x' \in    \{ z \in   { \mathbb R} ^N\,;\,  \vert  z-x \vert  =r,\,  \delta (z) \geq  \varepsilon _0r\,  \}$.

  \end{definition}

  In the above  ``quasi inwardly increasing near $\partial \Omega $"  means that $g$ is in some particular weak sense increasing along inward nontangential directions (near $\partial \Omega $) and "quasi" does not refer to a capacity.

  When $(QI)$ holds we also say that $g$ is $(C, \varepsilon _0)$ quasi inwardly increasing (near the boundary) in $\Omega $.

  Using property $(H)$ it is easy to see that -after a change of constant $C$- we may choose  the point $x'$ independently of $a>0$.

 It should also be noted that property ($QI$) is (up to change of constants)  a property of $g$ with respect  to  the pseudo-hyperbolic geometry in $\Omega $. To be more specific recall that for $c>0$, an arc $ \gamma :[0,\ell] \to  \Omega $ such that $ \delta ( \gamma (t) ) \geq c\, d( \gamma (t), \gamma (0))$ is called a $c$-John arc in $\Omega $. It can be shown that the validity of   ($QI$) (for some $C$ and $ \varepsilon _0$) means (up to change of constants) that  from  every point $x_0 \in   \Omega ( \varepsilon_0/4 )$ starts a $c$-John arc $ \gamma :[0,\ell) \to  \Omega $ in $\Omega $, such that $ \delta ( \gamma (\ell) )= \varepsilon _0$ and $ g(x,a)  \leq C g( \gamma (s), Ca)$ for all $a>0$ and $s \in   (0,\ell]$.

 \begin{theorem}  \label{thap2bis}  Suppose  that for some constants $C_2 \geq 1$ and $ \varepsilon _0>0$:

\quad  ($KOT$)$'$  $ \int _{2a}^ \infty  {\frac  { dt} {\sqrt{G(x,t)} } } \leq  C_2\, \sqrt{  {\frac  { a} {g(x,a)}} }$ whenever $a>0$, $x \in   \Omega $ and $g(x,a) \geq C_2\, a$,

\quad ($QI$) $g$ is $(C_2, \varepsilon _0)$-quasi inwardly increasing (near the boundary) in $\Omega $,

\noindent where $G$ is as in \eqref{cond1}. Then property (\ref{Potcond2bis}) holds for  some real $\overline c \geq 1$ depending only on $\Omega $ and $g$.\end{theorem}

 {\sl Remark.} For a more restricted class of $g$, a different and simpler approach to  this result  is  given  in {\bf \ref{sub4}}.

 Clearly ($KOT$)$'$ is -- up to change of constants -- a reinforcement of (KOT) in Proposition \ref{ConcCond-St}.   For the proof we require the following lemmas.

 \begin{lemma}  \label{lem2pr} Let $y_0 \in   \Omega $ and let $u$ be a positive solution of  (\ref{theequation}) in $\Omega $ that vanishes on $\partial \Omega  \cap B(y_0,2s)$,  $0<s \leq s_0$. Then, for   $z \in   \partial B(y_0,s) \cap \Omega $:
 $$u(z)\geq  { \frac {  1} {4 }} \sup  \{ u(y)\,;\,y \in   \Omega,\,  \vert  y-y_0 \vert  =s\, \} \quad \Rightarrow \quad d(z,\partial \Omega )  \geq   \eta s.$$
Here $s_0$ and $   \eta $ are (small) positive  constants depending only on  $\Omega $ and $g$. Moreover for all $z \in   \partial B(y_0,s)$ such that $ \delta (z) \geq   \eta  s$, we have, for a    constant $c'=c'(\Omega ,g) \geq 1$
$$u(z)\geq  { \frac {  1} {c' }} \sup  \{ u(y)\,;\,y \in   \Omega,\,  \vert  y-y_0 \vert  =s\, \}.$$
  \end{lemma}

{\sl Proof.}  The point here is that by Proposition \ref{ConcCond-St} and Remark \ref{schrrem} some uniform Harnack (or boundary Harnack) inequalities are available for $u$.

In particular, by the (inner uniform) Harnack inequalities satisfied by $u$, we may for the proof restrict to the case where  $B(y_0,2s) \cap    \partial \Omega \ne  \emptyset  $, i.e.\ $ \delta (y_0)<2s$.

Then, provided $s_0$ is chosen sufficiently small, and because  $\Omega $ is  bounded and Lipschitz,  there is a cylinder $T(r, \rho )$ in ${ \mathbb R} ^N$ -with $0<r=r(\Omega )< \rho= \rho (\Omega ) $, $ \rho /r=\kappa (\Omega )$- and an isometry $\mathfrak I$ of ${ \mathbb R} ^N$ such that $\omega (r, \rho ):= \mathfrak I(\Omega ) \cap T(r, \rho )$   is a standard Lipschitz domain (see Section 4) of height $ \rho $ and radius $r$ and $\mathfrak I (y_0)=(0, \tau ) \in   { \mathbb R} ^{N-1}\times { \mathbb R} \simeq { \mathbb R} ^N$ where $0< \tau \leq 4\kappa s$, $20\,s<r$.
Denote $A=\mathfrak I^{-1}(0,10\,  \kappa  s)$, $B=\mathfrak I^{-1}(0, \tau +s)$. By the Harnack boundary principle
$$u(x) \leq C\, u(A)\, s^{N-2} \mathbb G _A^V(x) \leq C\, u(A)\, s^{N-2} \mathbb G _A(x) \mathrm{ \ \ if \ }\mathfrak I( x)\in   \omega (5s, 5 \kappa  s)$$
where $\mathbb G _A^V$ is $L^V$ Green's function for $\mathfrak I^{-1}(\omega (10 s,10 \kappa s ))$, $\mathbb G _A:=\mathbb G _A^0$  and $C \geq 1$ is a constant depending only on $\Omega $ and $\bar c$. On the other hand by Harnack inequalities, $u(B) \geq c\, u(A)$ ($c \geq 0$ depending on $C_1$ and $\Omega $). Thus, when $z \in   \partial B(y_0,s) \cap    \Omega $
$$u(z) \leq C\, c^{-1} u(B)\, s^{N-2} \mathbb G _A(z)$$
But it is well-known that $s^{N-2} \mathbb G _A(z) \leq c' (d(z,\partial \Omega )/s)^ \alpha $ in $\mathfrak I^{-1}(\omega (5s, 5 { \frac {   \rho } {r }}s ))$ for some positive reals $c$ and $ \alpha $ depending only on $N$ and $ \rho /r$. The first part of the lemma follows.

The second part is then a simple consequence of the (uniform inner) Harnack inequalities that are satisfied by $L^V$ solutions. $\square$

\begin{lemma} \label{lemme3pr} Suppose $F \subsetneq   \partial \Omega $ is closed and $u=0$ in $\Omega \setminus F$. Given  $ \varepsilon _1>0$ there is a (large) $t=t( \varepsilon _1,\Omega ,g) \geq 2$ and a real $ \varepsilon _2>0$ such that $ \varepsilon _2t \leq {\frac  { 1} {10}}$ and

\hspace{5truemm} if $y_0 \in   \Omega $ and $ \delta (y_0) \leq  \varepsilon _2\,d(y_0,F)$ then $u(y_0) \leq  \varepsilon _1 \,\sup  \{ u(z)\,;\,  \vert  z-y_0 \vert  =t \delta (y_0)\, \}$.
\end{lemma}

{\sl Proof.} If $M=\sup  \{ u(z)\,;\,  \vert  z-y_0 \vert  =t \delta (y_0)\, \}$ the proof above shows, taking now $A= {\mathfrak I}^{-1}(0,\tau +2 t\delta (y_0)) $, that
$$u(y) \leq CM t^{-\alpha} $$ provided that $t \delta (y_0) \leq { \frac {  1} {10 }} d(y_0,F)$. The lemma follows. $\square$

 {\sl Proof of Theorem \ref{thap2bis}.} Let $u$ be a positive solution of (\ref{theequation}) and let $F$ be a compact subset of $\partial \Omega $ such that $u=0$ on $\partial \Omega \setminus F$. To prove the desired estimate we exploit again the ideas of the proof of Proposition \ref{ConcCond-St}. By (\ref{Potcond}) we may assume $F\ne  \emptyset  $ and $d(x,F)$ small.

 Consider a point $x_0 \in   \Omega $ and  the ball $B=B(x_0, R)$, $R:=d(x_0,F)< \varepsilon _0/10$.  Extending $u$ by zero outside $\Omega \setminus F$ we may view $u$ as a subharmonic function in $B$ and consider   the function
 $$w:=\sup  \{ u \circ { \mathcal R}\,;\, { \mathcal R}   \in   { \mathfrak I}_{ x_0}\, \}$$
 in $B$, where ${ \mathfrak I}_{x_0}$ is the group of isometries of $R^N$ fixing $x_0$.
Clearly $w$ is again  a continuous subharmonic function  in $B$  in the form $w(x)= \tilde u(d(x,x_0))$ for some continuous nondecreasing function $ \tilde u:[0, R) \to  { \mathbb R} _+$.

Denote, for $x \in   B$,
  ${ \mathcal U}(x):= \{  { \mathcal R}  \in   { \mathfrak I}_{x_0}\,;\,u({ \mathcal R}(x)) > { \frac {  1} {2 }} w(x) \}$. Note that $w(x)>0$.

Given $x_1 \in   B$, there is an open  neighborhood $\omega $ of $x_1$ such that, in $\omega $,     $u \circ { \mathcal R}<w$  for  $\mathcal R \in   { \mathfrak I}_{x_0} \setminus  { \mathcal U}(x_1)$, and  $u \circ { \mathcal R} \geq $ ${\frac  { 1} {4}} w$ if  ${ \mathcal R} \in   { \mathcal U}(x_1)$. So $w=\sup  \{ u \circ { \mathcal R}\,;\,  { \mathcal R} \in   { \mathcal U}(x_1)\,  \}$ in $\omega $.

 It follows that in $\omega $, the positive measure $ \Delta  w$ satisfies in the distribution sense
 \begin{align}
 \nonumber
 \Delta  w(x) & \geq \inf  \{  g({ \mathcal R}(x),u({ \mathcal R}(x)) \,;\, { \mathcal R} \in   { \mathcal U}(x_1) \}
 \\ \nonumber &\geq {\frac  { 1} {c_1}} \inf  \{ g(x_0,c_1^{-1} u({ \mathcal R}(x)))\,;\, { \mathcal R} \in   { \mathcal U}(x_1) \}
\\    \nonumber  &\geq  {\frac  { 1} {c_1}}  \  g(x_0,{\frac  { 1} {4c_1}}  w(x))
  \end{align}
where after the first inequality sign all terms are functions and where $c_1$ is a  new constant $ \geq 1$. For the second inequality, observe that for ${ \mathcal R} \in   { \mathcal U}_{x_1}$ and $x \in   \omega $,  $g(x_0,C ^{-1}u({ \mathcal R}(x))) \leq C g(x', u({ \mathcal R}(x)))$ for some  $x' \in   \partial B(x_0, \vert x-x_0 \vert  ) \cap \Omega $ such that $ \delta (x') \geq   \varepsilon _0 \vert  x-x_0 \vert $. So by Lemma \ref{lem2pr} and $(H)$, $g(x_0,C ^{-1}u({ \mathcal R}(x))) \leq c\,C g({ \mathcal R}(x), u({ \mathcal R}(x)))$ where $c=c(g,\Omega ) \geq 1$.

Thus,  $ \tilde u$ is a subsolution of the ODE  : $y''(r)+{ \frac {  N-1} {r }}y'(r)-g_0(y(r))=0$ in $(0, R)$ where $g_0(t):={ \frac {  1} {4c_1 }}g(x_0,{ \frac {  1} {4c_1 }}t)$.

Using Corollary \ref{Cor4} we get, if $g_0(u(x_0))>
0$ (and say $N \geq 3$), and if $G_0(t):=\int _0^tg_0(s)\, ds$, \medskip
 \begin{align}  \nonumber {\frac  { R} {N-2}}   {\frac  { 2^{N-2}-1} {2^{N-1}}  }   &  \leq \sqrt 2  \int _{u(x_0)}^{ \infty }  {\frac  { du} {\sqrt {G_0(u)-G_0(u(x_0))} } }
 \\ &=\sqrt {8c_1}\int _{u(x_0)}^{ \infty }  {\frac  { du} {\sqrt {G(x_0,u/4c_1)-G(x_0,u(x_0)/4c_1)} } } \nonumber
 \\ &=4\sqrt {8}\, (c_1)^{3/2}\,\int _{u(x_0)/4c_1}^{ \infty }  {\frac  { du} {\sqrt {G(x_0,u)-G(x_0,u(x_0)/4c_1)} } }.\nonumber
 \end{align}
 So if   $ g(x_0,u(x_0)/4c_1)>0$  it follows from (\ref{cond2}), ($KOT$)$'$ and the obvious adaptation of Lemma \ref{lemelem}  that
 $$c_4R \leq  \sqrt{  {\frac  { u(x_0)} {g(x_0,u(x_0)/4c_1)}} }$$
for some constant $c_4 >0$ depending only on $\Omega $ and $g$.

At this point we have established that for every $x_0 \in   \Omega$ and  some constant $\bar c=\bar c(\Omega ,g)$
$$g(x_0,u(x_0)/4c_1) \leq \bar c\, u(x_0)\,  d(x,F)^{-2}.$$

To get rid of the constant $4c_1$ in the r.h.s we use Corollary \ref{Cor4}  with  $\varepsilon _1=(4cCc_1)^{-1}$ where $c \geq 1$ is chosen as follows.  There is a (fixed) large $t$ and a small $   \varepsilon _2 >0$ (independent of $x_0$) with $ \varepsilon _2t \leq { \frac {  1} {10 }}$ and whenever  $ \delta (x_0) \leq  \varepsilon _2\,  d(x_0,F)$ then $u(x_0) \leq  \varepsilon _1\sup  \{ u(x)\,;\,  \vert  x-x_0 \vert  =t \delta (x_0)\,  \}$; so $u(x_0) \leq  c'\varepsilon _1\sup  \{ u(x)\,;\,  \vert  x-x_0 \vert  =t \delta (x_0),  \delta (x) \geq  \varepsilon _0t \delta (x)\,  \}$ (using property ($H$)). Here $c'=c'( \Omega , g, \varepsilon _0 ) \geq 1$. We choose $c:=c'$.

By the assumption ($QI$)  there exists $x_1 \in   \partial B(x_0,  t\delta (x_0))$ with $ \delta (x_1) \geq  \varepsilon _0t \delta (x_0)$ and $g(x_0,u(x_0)) \leq Cg(x_1, Cu(x_0))$. By the above, $u(x_0) \leq (4Cc_1)^{-1} u(x_1)$. Hence
 \begin{align}
   \nonumber g(x_0,u(x_0)) &\leq Cg(x_1,Cu(x_0)) \leq Cg(x_1,{ \frac {  1} {4c_1 }}u(x_1))
   \leq c {\frac  { u(x_1)} {  d(x_1,F)^2}}
\leq cc'' {\frac  { u(x_0)} { d(x_0,F)^2}}
  \end{align}
  where $c''$ is again a Harnack constant depending only on the chosen $t$, $\Omega $, $\varepsilon _0$ and $g$.

  This proves (\ref{Potcond2bis}) for $x \in   \Omega $ such that  $ \delta (x) \leq  \varepsilon d(x,F)$ if $ \varepsilon = \varepsilon (\Omega ,g)>0$ is small enough.  By (\ref{ConcCond-St}) the result then extends (with another constant)  to all $x \in   \Omega $.
\qed

\medskip

\subsection  {\large\bf Few Examples. } \label{ex3}

We first describe some examples of functions $g$ in the form $g(x,t)=b(x)\, \overline g(t)$ that satisfy the basic assumptions in section 1 (we then say that $g$ is ``good"), in particular  \eqref{Vcond} by the results above. We owe to Moshe Marcus the idea of considering examples of this form. See also the reference  \cite{MVcpam}.

a) Recall that $\Omega $ always denote a bounded Lipschitz domain in ${ \mathbb  R}^N$. Suppose  that $\overline g$ is convex and nondecreasing in $[0, \infty )$ with $g(0)=0$ and  that $b \in   C_+(\Omega )$. Suppose furthermore that   for some constants $ \varepsilon >0$ and $c \geq 1$,

(i) if $x,\, x' \in   \partial \Omega $, $ \vert  x-x' \vert   \leq { \frac {  1} {2 }}\, d(x,\partial \Omega) $ then $b(x') \leq cb(x)$,

(ii) whenever  $ \zeta   \in   \partial \Omega $, $0<2r <r' \leq  \varepsilon $, $x \in   \partial B( \zeta  ,r) \cap \Omega $,  there exists a point $x' \in   \partial B( \zeta  ,r')$ such that $b(x) \leq c\, b(x')$ and $d(x',\partial \Omega ) \geq  \varepsilon \, r'$. \

(iii)  $\int _a^ \infty { \frac {  dt} {\sqrt{\overline G(t)} }} \leq c\, \sqrt{{ \frac {  a} {\overline g(a) }}}$ for $a$ large -where $\overline G(t):=\int _0^t\, \overline g(t)\, dt$.

Then $g(x,t)$ is ``good" and satisfies the assumptions of Theorem \ref{thap2bis}. Note that by the main results in  \cite{Annals87} any positive $ \beta  \in   C_+(\Omega )$ which  solves $ \Delta  u-Vu=0$ near $\partial \Omega $, with $ \sup _\Omega  \delta ^2 \vert  V  \vert < \infty $, $V_-=o( \delta ^2) $ as $ \delta  \to  0$, and $ b =0$ on $\partial \Omega $ is an admissible function $b$ as well as any power $ \beta ^m$ of $ \beta $, $m>0$. Clearly more general second order elliptic operators can be considered.

 To obtain more examples notice that $b(x)= \delta (x) $ is an admissible function $b$, that powers, products, sums of admissible functions $b$ are still admissible.

{\sl Example 1.}  $g(x,t)= c\,\varphi  (x)^ \alpha\,  \delta (x)^ \beta t^q$ where $ \varphi  $ is the ground state of the Dirichlet Laplacian in $\Omega $ and  $ \alpha,\,  \beta   \geq 0$.

{\sl Example 2.} It is easily checked that any increasing function $\bar g:{ \mathbb  R}_+ \to  { \mathbb  R}_+$ such that $g(t) \geq c ({ \frac {  t} {s }})^ \alpha g(s)$ for some $ \alpha >1$, a constant $c>0$ and all $t>s>1$ is an admissible function $\bar g$.

\medskip

b) Other examples can be obtained using the following remark: let  $g_1,\dots, \, g_n$ be admissible functions $g$ in $\Omega \times { \mathbb  R}_+$, in the sense that each satisfies  the assumptions of Proposition \ref{ConcCond-St} (resp. Theorem \ref{thap2bis}) Then if $ \beta _1,\dots, \,  \beta _n$ are positive reals, the function $ \sum_{j \leq N}  \beta _jg_j$ is  again an admissible function $g$. The proof is left as an exercise.

c) We finally mention other classes of examples. If $g(x,t):=\overline  g (  t\,b_1 (x)) \,b_2 (x)$ where $\overline g:{ \mathbb  R}_+ \to  { \mathbb  R}_+$, $b_1,\,b_2 :\Omega  \to  { \mathbb  R}_+$ are admissible in the sense of a) above, then, provided $\overline g$ satisfies the $ \Delta  _2$ condition, $g$ satisfies the assumptions of Theorem \ref{thap2bis}.

 \subsection{\label{sub4}} \label{A44}  {\large\bf A variant of Theorem \ref{thap2bis} }

 In this part, we describe another way to obtain estimates in the form (\ref{Potcond2bis} ).  We show that under condition$(QI)$ and another mild condition on $g$, property (\ref{Potcond}) implies property (\ref{Potcond2bis}) for solutions.

 \begin{theorem}  \label{thap2} Consider  in a bounded Lipschitz domain $\Omega $ of  ${ \mathbb R} ^N$ the equation
\begin{equation} \label{same}\Delta  u(x) =g(x,u(x))\end{equation}
where   $g:\Omega \times { \mathbb  R}_+ \to  { \mathbb  R}_+$ is continuous and such that $t \mapsto h(x,t)=t^{-1}g(x,t)$ is nondecreasing in ${ \mathbb  R}_+^*$.
Suppose  that for some constants $C \geq 1$ and $ \varepsilon _0>0$

\quad (i) ($QI$) $g$ is $(C, \varepsilon _0)$-quasi inwardly increasing (near the boundary) in $\Omega $.

\quad (ii) every positive solution $u$ of (\ref{same})
is such that  $g(x,u(x)) \leq C\, ( \delta (x))^{-2}\, u(x)$ in $\Omega $.

 Then property (\ref{Potcond2bis}) holds for $\Omega $ and the equation (\ref{same}).\end{theorem}

\vspace{2mm}

 {\sl Proof.} a) Let $u$ be a positive solution of (\ref{same}) in $\Omega $ and let $F$ be the smallest closed subset of $\partial \Omega $ such that $ u=0$ in $\partial \Omega \setminus F$.
 By the assumption we have that $u$ is a solution of an equation $ \Delta  u-Vu=0$ in $\Omega $ where the potential $V \in  L^ \infty _{loc}(\Omega )$ is nonnegative and such that
 $$V(x) \leq { \frac { C} {  \delta (x)^2}} \mathrm{ \ for\ } x  \in   \Omega. $$

 b) Suppose that for some cylinder $T(r, \rho )= \{ (x',x_N) \in   { \mathbb  R}^{N-1}\times { \mathbb  R}\simeq { \mathbb  R}^N\,;\,  \vert  x' \vert  <r,\;  \vert  x_N \vert  < \rho \, \}$ with $0<r< \rho $, the region $\omega (r, \rho ):=\Omega  \cap T(r, \rho )$ is a standard Lipschitz domain of radius $r$ and height $ \rho $. Suppose also that $u=0$ on $\partial \Omega   \cap T(r, \rho )$ (i.e. $F \cap T(r, \rho )= \emptyset  $) and denote $A'=(0,{ \frac {  3} {4 }}\, \rho )$.

 Then by the Harnack boundary property,
 $$u(x) \leq c\, r^{N-2}\,u(A')\,\mathbb G ^V_{A'}(x) \leq c\, r^{N-2}\,u(A')\,\mathbb G _{A'}(x)$$ for $x \in   \Omega  \cap T(r/2, \rho /2)$ and some $c=c(C,{\frac  { r} { \rho }}, N) \geq 1$. Here $\mathbb G ^V$ is the $L^V$ Green's function in $\omega :=T(r, \rho ) \cap \Omega $, $\mathbb G =\mathbb G ^0$ and $A'=(0,{ \frac {  3} {4 }} \rho )$.

 Thus on using the Harnack inequalities for $L^V$ we get for $x' \in   \partial B(x,r/2)$, $ \delta (x') \geq  \varepsilon _0r/2$,
 $$u(x) \leq  c\,c'\, r^{N-2}\,u(x')\,\mathbb G _{A'}(x)$$ with a constant $c'$ depending only on $C$, $N$ and $ \rho /r$.

Since $r^{N-2} \mathbb G _{A'}(x) \leq c_1\, (d(x,\partial \Omega )/r)^ \alpha $ when  $x \in   \Omega  \cap T(r/2, \rho /2)$ for some constants $c_1 \geq 1$ and $ \alpha  \in   (0,1)$, it follows that $$u(x) \leq     { \frac {  1} {C }}\,u(x')$$ if  in addition to the previous conditions $x \in   T(  \varepsilon _1r, \varepsilon _1\, \rho ) \cap \Omega $ where the  constant $ \varepsilon _1>0$ is chosen small enough depending only on $c$, $c'$, $N$ and ${ \frac {   \rho } {r }} $.

c) The domain  $\Omega $ being bounded and  Lipschitz, there is a small $\kappa >0$ such that for every point $x \in   \Omega $ for which  $ \delta (x) \leq \kappa d(x,F)$ the following holds: $d(x,\partial \Omega ) \leq { \frac {  1} {8 }} \varepsilon _1 \varepsilon _0$ and there is an affine isometry ${ \mathcal R}$ of ${ \mathbb  R}^N$ and numbers $r$, $ \rho >0$ such that ${ \mathcal R}^{-1} (\Omega )  \cap T(r, \rho )$ is a standard Lipschitz domain of radius $r$ and  height $ \rho $, $F \cap   { \mathcal R}(T(r, \rho ))= \emptyset   $, the point $x $ being inside   $ { \mathcal R}( T(\varepsilon _1r,\varepsilon _1 \rho ))$. Moreover the ratio ${ \frac {   \rho } {r }}$ depends only on $\Omega $ and one may further choose $r \geq c_2\, d(x, F)$ for some $c_2>0$.

 Thus if $x' \in   \Omega \cap \partial B(x,r/2) $ and $ \delta (x') \geq  \varepsilon _0r/2$, we have by the assumptions $g(x,a) \leq Cg(x',Ca)$ for all $a>0$ -or $h(x,a) \leq C^2h(x',Ca)$- and
 $$V(x)=h(x,u(x)) \leq h(x, { \frac {  1} {C }}u(x'))  \leq C^2\, h(x',  u(x'))  \leq { \frac {2\,C^3 } { r^2}} \leq \, { \frac { 2\, C^3 (c_2)^{-2} } { d(x,F)^2}}.$$
 Here we have used the fact that $h(x',  u(x'))   \leq C/( \delta (x'))^2 \leq 2\,C/r^2$.

On the other hand if $x \in   \Omega $ is such that $\delta (x) \geq \kappa d(x,F)$  then $$V(x) \leq { \frac {  C} {  \delta (x)^2}} \leq  {\frac { C\, \kappa^{-2}} {  d(x,F)^2}}$$

 and the proposition follows. $\square$

\end{document}